\documentclass[11pt,bezier]{article}
\usepackage{amsmath,amssymb,amsfonts,euscript,graphicx,rotating}
\usepackage[usenames]{color}
\usepackage[all]{xy}
\newpage

\newcommand{\be}{\begin{equation}}
\newcommand{\ee}{\end{equation}}
\newtheorem{theorem}{Theorem}[section]

\newtheorem{lemma}{Lemma}[section]

\newtheorem{example}{Example}[section]

\newtheorem{remark}{Remark}[section]

\title{\bf\Large Two New Estimators of Entropy for Testing Normality
}
\vspace{-1cm}
\author
{Akram Kohansal, Saeid Rezakhah \footnote{Faculty of Mathematics and Computer Science, Amirkabir University of Technology,
Tehran, Iran. Email: ak\_kohansal@aut.ac.ir, rezakhah@aut.ac.ir}}
\begin{document}
\maketitle
\begin{abstract}
\par
We present two new estimators for estimating the entropy of absolutely continuous random variables. Some properties of them are considered, specifically consistency of the first is proved. The introduced estimators are compared with the existing entropy estimators. 
Also, we propose two new tests for normality based on the introduced entropy estimators and compare their powers with the powers of other tests for normality. The results show that the proposed estimators and test statistics perform very well in estimating entropy and testing normality. A real example is presented and analyzed.\\ \quad \\
{\it Keywords}: Information theory, Entropy estimator, Moving average method, Ranked set sampling, Testing normality. \\ \quad \\
{\it Mathematics Subject Classification:} 62G10, 62G30.

\end{abstract}

\section{Introduction}
\par
Suppose a random variable $X$ has a distribution function $F(x)$ with an absolutely continuous density function $f(x)$. The entropy $H(f)$ of the random variable is defined by Shannon \cite{30me} to be
\begin{equation}\label{H}
H(f)=\int_{-\infty}^{\infty}f(x)\log f(x)\mathrm{d}x.
\end{equation}
The non-parametric estimation of $H(f)$ has been discussed by many authors, including Vasicek \cite{32me}, Van Es \cite{31me}, Ebrahimi {\em et al.} \cite{12me}, Correa \cite{8me} and Wieczorkowski and Grzegorewski \cite{33me}.
\par
Among these various entropy estimators, Vasicek's sample entropy has been most widely used in developing entropy-based goodness-of-fit tests, including Ebrahimi {\em et al.} \cite{13me}, Park and Park \cite{25me} and Alizadeh Noughabi \cite{3me} .
\par
Vasicek's estimator was based on the fact that the equation (\ref{H}) can be expressed as
\begin{equation}\label{F}
H(f)=\int_{0}^{1}\log\left\{\frac{\mathrm{d}}{\mathrm{d}p}F^{-1}(p)\right\}\mathrm{d}p.
\end{equation}
The estimate was obtained by replacing the distribution function $F$ by the empirical distribution function $F_n$, and using a difference operator instead of the differential operator. The derivative of $F^{-1}(p)$ is then estimated by a function of the order statistics.
\par
Let $X_1,\dots,X_n,~n\geq3$ be a sample from the distribution $F$. Let $X_{(1)}\leq\dots\leq X_{(n)}$ denote the order statistics from the sample $X_1,\dots,X_n$. The Vasicek's estimator of entropy has a following form
\begin{equation*}
HV_{mn}=\frac{1}{n}\sum_{i=1}^{n}\log\left\{\frac{n}{2m}\left(X_{(i+m)}-X_{(i-m)}\right)\right\},
\end{equation*}
where the window size $m$ is a positive integer smaller than $n/2$, $X_{(i)}=X_{(1)}$ if $i<1$, $X_{(i)}=X_{(n)}$ if $i>n$. Vasicek proved that his estimator is consistent, i.e. ${HV_{mn} \stackrel{Pr.}{\longrightarrow} H(f)}$ as $n,m\rightarrow\infty,~\frac{m}{n}\rightarrow0$.
\par
Van Es \cite{31me} proposed another estimator of entropy and established the consistency and asymptotic normality of this estimator under some conditions. Van Es's estimator is given by
\begin{equation}\label{VE}
HVE_{mn}=\frac{1}{n-m}\sum_{i=1}^{n-m}\log\left\{\frac{n+1}{m}\left(X_{(i+m)}-X_{(i)}\right)\right\}+\sum_{k=m}^{n}\frac{1}{k}+\log(m)-\log(n+1).
\end{equation}
\par
Ebrahimi {\em et al.} \cite{12me} adjusted the weights of Vasicek's estimator, in order to take into account the fact that the differences are truncated around the smallest and the largest data points (i.e. $X_{(i+m)}-X_{(i-m)}$ is replaced by $X_{(i+m)}-X_{(1)}$ when $i\leq m$ and $X_{(i+m)}-X_{(i-m)}$ is replaced by $X_{(n)}-X_{(i-m)}$ when $i\geq n-m+1$). Their estimator is given by
\begin{equation*}
HE_{mn}=\frac{1}{n}\sum_{i=1}^{n}\log\left\{\frac{n}{c_im}\left(X_{(i+m)}-X_{(i-m)}\right)\right\},
\end{equation*}
where
\begin{equation*}
c_i=\left\{\begin{array}{cc}
                                                                         1+\frac{i-1}{m} & 1\leq i\leq m, \\
                                                                         2               & m+1\leq i\leq n-m, \\
                                                                         1+\frac{n-i}{m} & n-m+1\leq i\leq n. \\
\end{array}\right.
\end{equation*}
They proved that ${HE_{mn} \stackrel{Pr.}{\longrightarrow} H(f)}$ as $n,m\rightarrow\infty,~\frac{m}{n}\rightarrow0$. They compared their estimator with Vasicek's estimator and by simulation, showed that their estimator has a smaller bias and mean squared error.
\par
Correa \cite{8me} proposed a modification of Vasicek's estimator. In estimation the density $f$ of $F$ in the interval $\left(X_{(i-m)},X_{(i+m)}\right)$ he used a local linear model based on $2m+1$ points: $F(X_{(j)})=\alpha+\beta X_{(j)}+\varepsilon,~j=m-i,\dots,m+i$.
This yields a following estimator
\begin{equation*}
HC_{mn}=-\frac{1}{n}\sum_{i=1}^{n}\log\left(a_i\right),
\end{equation*}
where $\displaystyle a_i=\frac{\sum_{j=i-m}^{i+m}(X_{(j)}-\bar{X}_{(i)})(j-i)}{n\sum_{j=i-m}^{i+m}(X_{(j)}-\bar{X}_{(i)})^2},$ and
\begin{equation}\label{xbar}
\bar{X}_{(i)}=\frac{1}{2m+1}\sum_{j=i-m}^{i+m}X_{(j)}.
\end{equation}
He compared his estimator with Vasicek's estimator and Van Es's estimator. The MSE of his estimator is consistently smaller than the MSE of Vasicek's estimator. Also, for some of $m$, his estimator behaves better than the Van Es's estimator.
\par
Wieczorkowski and Grzegorzewsky \cite{33me} provided a modification of Vasicek estimator and Correa estimator. Their estimator is given by
\begin{equation}\label{W}
HW_{mn}=HV_{mn}-\log(n)+\log(2m)+c,
\end{equation}
where
\begin{equation}\label{c}
c=-\left(1-\frac{2m}{n}\right)\Psi(2m)+\Psi(n+1)-\frac{2}{n}\sum_{i=1}^{m}\Psi(i+m-1),
\end{equation}
and $\Psi(x)$ is the digamma function defined by $\Psi(x)=\frac{\mathrm{d}\log\Gamma(x)}{\mathrm{d}x}=\frac{\Gamma'(x)}{\Gamma(x)},$ and for integer arguments
$\Psi(k)=\sum_{i=1}^{k-1}\frac{1}{i}-\gamma,$ where $\gamma$ is the Euler constant, $\gamma=0.57721566490\dots$.
\par
Zamanzadeh and Arghami \cite{36me} provided two new entropy estimators as follows:
\begin{eqnarray}
HZ1_{mn}=\frac{1}{n}\sum_{i=1}^{n}\log(b_i),\nonumber\quad
HZ2_{mn}=\sum_{i=1}^{n}w_i\log(b_i),\nonumber
\end{eqnarray}
where
\begin{equation}\label{bi}
b_i=\frac{X_{(i+m)}-X_{(i-m)}}{\sum_{j=k_1(i)}^{k_2(i)-1}\left(\frac{\hat{f}(X_{j+1})+\hat{f}(X_{j})}{2}\right)\left(X_{(j+1)}-X_{(j)}\right)},
\end{equation}
and $k_1(i)$ is $1$ and $i-m$ for $i\leq m$ and $i>m$, respectively. Also, $k_2(i)$ is $i+m$ and $n$ for $i\leq n-m$ and $i>n-m$, respectively. Further, $\hat{f}(X_i)=\frac{1}{nh}\sum_{j=1}^nk\left(\frac{X_i-X_j}{h}\right),$ where $k$ is chosen to be Normal density function and $h=1.06\hat{\sigma}n^{-1/5},$ in which $\hat{\sigma}=\{\frac{1}{n}\sum_{i=1}^n(X_i-\bar{X})^2\}^{1/2},$ and
\begin{equation}\label{wi}
w_i=\left\{\begin{array}{ccc}
  \frac{m+i-1}{\sum_{i=1}^nw_i} & 1\leq i\leq m, &  \\
  \frac{2m}{\sum_{i=1}^nw_i} & m+1\leq i\leq n-m, & i=1,\dots,n. \\
  \frac{n-i+m}{\sum_{i=1}^nw_i} & n-m+1\leq i\leq n, &
\end{array}
\right.
\end{equation}

\par
Ranked Set Sampling (RSS) has been under vast investigations since its introduction. Patil {\em et al.} \cite{2ma} comprehensively reviewed works done on RSS. Takahasi and Wakimoto \cite{3ma} and Dell and Clutter \cite{4ma} were the first authors who examined the theory of RSS mathematically. Stokes \cite{5ma} proposed an RSS estimator of the population variance and showed its asymptotic unbiasedness regardless of the presence of errors in the ranking. MacEachern {\em et al.} \cite{6ma} developed an unbiased estimator of the variance of the population based on an RSS and demonstrated that this new estimator is more efficient than that of Stokes. Stokes and Sager \cite{7ma} studied the properties of the empirical distribution function based on RSS and showed that it is unbiased and has greater precision than that of Simple Random Sampling (SRS). {\textcolor{blue} { Balakrishnan and Li  \cite{bala_ORSS} proposed ordered ranked set sampling (ORSS) and
 developed optimal linear inference based on ORSS. Recently, Balakrishnan {\em et al.} \cite{bala2} considered the goodness-of-fit procedure for testing  normality based on the empirical characteristic function for ranked set sampling data. A very well written and complete overview of the available results in the context of RSS is the book of Chen {\em et al.} \cite{balabook}, which includes 175 key references on this subject. }}

\par
The objective of this article is two fold. First, we develop the entropy estimators in order to obtain the estimators with less bias and less mean squared error. They are also more efficient than others. These estimators are obtained by using the moving average and ranked set sampling methods; and modifying the Van Es's \cite{31me} estimator and Wieczorkowski and Grzegorzewsky's \cite{33me} estimator. We also consider the scale invariance property of variances and mean squared errors of them. Moreover, the consistency of the first estimator is established. The second objective of this article is to improve the tests for normality in order to obtain the most powerful tests. These test statistics are provided based on the proposed entropy estimators. We also compare the power of these tests with the other tests, using Monte Carlo computations. Differences in power of the tests are considerable, but the results show that the first test is most powerful against the alternatives with support $(0,\infty)$ such as Gamma, Weibull, Log normal, Beta and the second is most powerful against the alternatives with support $(-\infty,\infty)$ such as t, Extreme value, Logistic.

\par
The rest of this paper is arranged as follows: In Section 2, we first enhance the moving average and ranked set sampling methods; and introduce two new entropy estimators. Some properties of them are studied in the same section. Also, we compare our estimators with the competitor estimators by a simulation study.  In Section 3, we introduce goodness-of-fit tests for normality, based on the proposed entropy estimators and then compare their powers with the powers of other tests of normality. A real example is also presented and analyzed in Section 4.


\section{New entropy estimators}
\par
In this section, we introduce two entropy estimators and compare them with other estimators.

\subsection{Moving average method}
\par
In statistics, smoothing a data set is to create an approximating function that attempts to capture important patterns in the data, while leaving out noise phenomena. One of the most common smoothing methods is moving average. This method is a technique that can be applied to the time series analysis, either to produce smoothed periodogram of data, or to make forecasts \cite{11me}.
\par
A moving average (MA) method is the unweighted mean of the previous $n$ datum points. Suppose individual observations, $X_1,\dots,X_n$ are collected. The moving average of width $w$ at time $i$ is defined by
$$Y_i=\frac{X_i+X_{i-1}+\dots+X_{i-w+1}}{w}=\frac{\sum_{j=i-w+1}^iX_{j}}{w}\quad i\geq w.$$
For periods $i<w$, we do not have $w$ observations to calculate a moving average of width $w$.
\par
Now, we develop the construction of the MA method. For this aim, we defined the moving average of width $w$, which is an odd integer, at time $i$ as:
\begin{equation*}\label{rev}
Y_i=\frac{X_{i-\frac{w-1}{2}}+X_{i-\frac{w-3}{2}}+\dots+X_{i+\frac{w-1}{2}}}{w}~~~~{(w+1)}/{2}\leq i\leq n-{(w-1)}/{2}.
\end{equation*}
For $i \geq n-\frac{w-3}{2}$ and $i \leq \frac{w-1}{2}$, the moving average at time $i$ is defined as the average of all observations that are equal or greater than $X_{i}$ and equal or smaller than $X_{i}$, respectively. So, we define $Y_i$ as follows:

\begin{equation}\label{Y}
Y_i=\left\{
\begin{array}{lll}
  \frac{\sum_{j=1}^{i}X_{j}}{i} & 1\leq i \leq {(w-1)}/{2}, \\
  \frac{\sum_{j=i-\frac{w-1}{2}}^{i+\frac{w-1}{2}}X_{j}}{w} & {(w+1)}/{2}\leq i\leq n-{(w-1)}/{2}, & i=1,\dots,n.\\
  \frac{\sum_{j=i}^{n}X_{j}}{n-i+1} & n-{(w-3)}/{2}\leq i \leq n,
\end{array}
\right.
\end{equation}

\par
One characteristic of the MA is that if the data have an uneven path, applying the MA will eliminate abrupt variation and cause the smooth path.

\subsection{Ranked set sampling method}
\par
RSS is known to be a statistical method for data collection that generally leads to more efficient estimators than competitors based on SRS. The concept of RSS was used first time by McIntyre \cite{6do}, to estimate the population mean of pasture yields in agricultural experimentation. Since then there has been substantial progress in studying this sampling scheme and its extensions. When the variable of interest can be more readily ranked than measured, RSS provides improved statistical inference. We now briefly explain the concept of RSS for completeness.
\begin{itemize}
  \item [$\bullet$] $n$ random samples, each of size $n$, are drawn from the population.
  \item [$\bullet$] The $i$th sample ($i=1,\dots,n$) is inspected to identify the unit of (judgement) $i$th lowest rank.
  \item [$\bullet$] Finally, the $n$ identified units in the (judgement) ranked set are measured.
\end{itemize}

\par
In the next subsection, we use the MA and RSS methods; and present two new entropy estimators are presented.

\subsection{Entropy estimators}
\par
Suppose $X_1,\dots,X_n$ are an SRS from an unknown absolutely continuous distribution $F$ with a probability density function $f(x)$. The simulation study, which provided with some papers such as Wieczorkowski and Grzegorzewsky \cite{33me} and Alizadeh Noughabi \cite{3me}, shows that in most cases $HVE_{mn}$, which introduced in (\ref{VE}), and $HW_{mn}$, which introduced in (\ref{W}), have the least Root Mean Squared Errors (RMSEs) and Standard Deviations (SDs) in estimating entropy.
\par
Now, we use the MA and RSS methods to modify $HVE_{mn}$ and $HW_{mn}$; and obtain two new entropy estimators. According to (\ref{F}), we know
\begin{equation*}
H(f)=\int_{0}^{1}\log\left\{\frac{\mathrm{d}}{\mathrm{d}p}F^{-1}(p)\right\}\mathrm{d}p,
\end{equation*}
$F^{-1}(p)$ as a function of quantiles in previous equation is the sample path of order statistics, but usually it is not smooth. So we propose to imply the MA method of proper order, say $w$, to smooth this sample path.
\par
To obtain the new estimators, we select $n$ SRS and repeat $n$ times this work and symbolize the SRS in $i$th times by $X_1^i,\dots,X_n^i$. Also, $X_{i1}<\dots<X_{in}$ denote the order statistics from the sample $X_1^i,\dots,X_n^i$. By using the MA method in each times, we define the new variables from (\ref{Y}) as
\begin{equation}\label{YY}
Y_{ki}=\left\{
\begin{array}{lll}
  \frac{\sum_{j=1}^{i}X_{kj}}{i} & 1\leq i \leq {(w-1)}/{2}, \\
  \frac{\sum_{j=i-\frac{w-1}{2}}^{i+\frac{w-1}{2}}X_{kj}}{w} & {(w+1)}/{2}\leq i\leq n-{(w-1)}/{2}, & k,i=1,\dots,n.\\
  \frac{\sum_{j=i}^{n}X_{kj}}{n-i+1} & n-{(w-3)}/{2}\leq i \leq n,
\end{array}
\right.
\end{equation}
 {\textcolor{blue} { Following  the method of Order Ranked Set Sampling (ORSS), \cite{bala_ORSS}, we choose $Y_{11},\dots,Y_{nn}$ and symbolize the order statistics of them as $Y^R_1,\dots,Y^R_n$.}} So, the new estimators are defined as:
\begin{eqnarray}
\label{VER}HVE^{R}_{mn}&=&\frac{1}{n-m}\sum_{i=1}^{n-m}\log\left(Y^R_{i+m}-Y^R_{i}\right)+\sum_{k=m}^{n}\frac{1}{k},\\
\label{WR}HW^{R}_{mn}&=&\frac{1}{n}\sum_{i=1}^{n}\log\left(Y^R_{i+m}-Y^R_{i-m}\right)+c,
\end{eqnarray}
where $c$ is defined in (\ref{c}), $m$ is a positive integer smaller than $n/2$ and $Y^R_i=Y^R_1$ if $i<1$, $Y^R_i=Y^R_n$ if $i>n$.
\par
We can easily prove that the scale of the random variable has no effect on the accuracy of $HVE^{R}_{mn}$ and $HW^{R}_{mn}$ in estimating $H(f)$.
\begin{remark}
Let $H^{W_1}(f)$ and $H^{W_2}(f)$ denote entropies of the distribution of continuous random variables $W_1$ and $W_2$, respectively, and $W_2=kW_1$, where $k>0$. It is easy to see that $\displaystyle H^{W_2}_i=\log k+H^{W_1}_i,~ i=1,2.$ Then the followings hold
\begin{itemize}
\item {$E(H^{W_2}_i)= E(H^{W_1}_i)+\log k,~i=1,2,$}
\item {$Var(H^{W_2}_i)= Var(H^{W_1}_i),~i=1,2,$}
\item {$MSE(H^{W_2}_i)= MSE(H^{W_1}_i),~i=1,2,$}
\end{itemize}
where $H^{W_i}_1=HVE^R_{mn}$ and $H^{W_i}_2=HW^R_{mn}$ in which the superscript $W_i,i=1,2$ refer to the corresponding distribution. $\blacksquare$
\end{remark}

\par
\begin{example}{\em
For the explanation of the MA method, we simulate 30 samples from the standard Normal distribution and plot their order statistics in Figure $1$ with $A$.
The sample path of order statistics is smoothed by MA of order $3$. New variables are defined from (\ref{YY}) and the smoothed path of new variables is plotted in Figure 1 with $B$. This plot shows that the new sample path is smoother than the sample path of the original order statistics. Also, with considering MA of order 5, we define new variables from (\ref{YY}) and plot them in Figure 1 with $C$. Even though the smoothing sample path of order statistics by using the MA of order 3 is not as smooth as using MA of order 5, the plot $B$ is very similar to real plot i.e. $A$. So without loss of generality, we just consider MA of order $w=3$ in (\ref{YY}). $\blacksquare$}
\end{example}

\input{epsf}
\epsfxsize=5in \epsfysize=4in
\begin{figure}
\label{order}
\centerline{\epsfxsize=4in \epsfysize=3in \epsffile{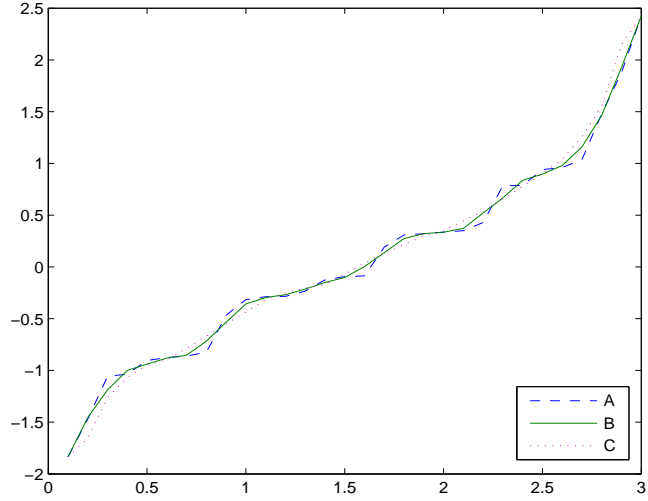}}
\vspace{-0.3in}
\caption{\footnotesize Sample path of order statistics from the standard Normal distribution (A), smoothed path of order 3 (B) and order 5 (C) at MA method. }
\end{figure}


  \begin{lemma}
{\textcolor{blue} {Let $C$ be the class of continuous densities with finite entropies, and let $X_1,\dots,X_n$ be a random sample with some density $f\in C.$   Then $HVE^{R}_{mn}-HVE_{mn}\rightarrow0$} for all $m=o(n)$.}
\end{lemma}
{\bf\em {Proof:}} {\textcolor{blue} { We show that $\left|HVE^{R}_{mn}-HVE_{mn}\right|\rightarrow 0$ then  Squeeze theorem implies that $HVE^{R}_{mn}-HVE_{mn}\rightarrow 0$. To visualize such convergence we use the following inequality 
\begin{align*}
0\leq\left|HVE^{R}_{mn}-HVE_{mn}\right|
=&\bigg|\frac{1}{n-m}\sum_{i=1}^{n-m}\log(Y^R_{i+m}-Y^R_{i})-\sum_{i=1}^{n-m}\log(X_{(i+m)}-X_{(i)})\bigg|\\
=&\bigg|\underbrace{\frac{1}{n-m}\sum_{i=1}^{n-m}\log\left(\frac{Y^R_{i+m}-Y^R_{i}}{X_{(i+m)}-X_{(i)}}\right)}_{A(n,m)}\bigg|,
\end{align*}
where while $A(n,m)>0$,
\begin{align*}
\left|A(n,m)\right|=&\frac{1}{n-m}\sum_{i=1}^{n-m}\log\left(\frac{Y^R_{i+m}-Y^R_{i}}{X_{(i+m)}-X_{(i)}}\right)
\leq\frac{1}{n-m}\sum_{i=1}^{n-m}\log\left(\frac{X_{i+m+1,i+m+1}-X_{i-1,i-1}}{X_{(i+m)}-X_{(i)}}\right)\\
\leq& \frac{1}{n-m}\sum_{i=1}^{n-m}\log\left(\frac{X_{i+m+1,i+m+1}-X_{i-1,i-1}}{X_{(i+m-k)}-X_{(i-1)}}\right)\;\;\;\mbox{for some $k \in {\mathbb N}$}\\
=&\frac{1}{n-m}\sum_{i=1}^{n-m}\log\left(\frac{X_{i+m+1,i+m+1}-X_{i-1,i-1}}{F(X_{i+m+1,i+m+1})-F(X_{i-1,i-1})}\right.\\
&\qquad\qquad\left.\times \frac{F(X_{(i+m-k)})-F(X_{(i-1)})}{X_{(i+m-k)}-X_{(i-1)}}\times \frac{F(X_{i+m+1,i+m+1})-F(X_{i-1,i-1})}{F(X_{(i+m-k)})-F(X_{(i-1)})}\right)\\
\longrightarrow&\frac{1}{n-m}\sum_{i=1}^{n-m}\log\left(\frac{1}{f(X_{j'_{i}})}\times f(X_{j_{i}})\times1\right)\longrightarrow0.
\end{align*}
The second inequality is valid since there exists $k>0$ such that $X_{(i+m-k)}-X_{(i-1)}\leq X_{(i+m)}-X_{(i)}$ and  $m+i-k\rightarrow\infty$ as $m\rightarrow\infty$. Also, the first convergence is valid because there exist a value $X_{j_{i}}\in (X_{(i-1)},X_{(m')})$ as $m'\rightarrow\infty$ such that $\frac{F(X_{(m')})-F(X_{(i-1)})}{X_{(m')}-X_{(i-1)}}=f(X_{j_{i}})$. Moreover, the last convergence is valid by the by the Ces\`{a}ro mean theorem.
If $A(n,m)<0$ then
\begin{align*}
\left|A(n,m)\right|=&\frac{1}{n-m}\sum_{i=1}^{n-m}\log\left(\frac{X_{(i+m)}-X_{(i)}}{Y^R_{i+m}-Y^R_{i}}\right)
\leq\frac{1}{n-m}\sum_{i=1}^{n-m}\log\left(\frac{X_{(i+m)}-X_{(i)}}{X_{i+m-1,i+m-1}-X_{i+1,i+1}}\right)\\
\leq& \frac{1}{n-m}\sum_{i=1}^{n-m}\log\left(\frac{X_{(i+m+k)}-X_{(i+1)}}{X_{i+m-1,i+m-1}-X_{i+1,i+1}}\right)\;\;\;\mbox{for some $k \in {\mathbb N}$}.\\
\end{align*}
The last inequality is valid since there exists $k>0$ such that $X_{(i+m+k)}-X_{(i+1)}\leq X_{(i+m)}-X_{(i)}$ and  $m+i+k\rightarrow\infty$ as $m\rightarrow\infty$. The remaining of the proof is similar to the previous case. $\blacksquare$}}

\begin{theorem}
{\textcolor{blue} {Let $C$ be the class of continuous densities with finite entropies, and let $X_1,\dots,X_n$ be a random sample from $f\in C.$ If $m=o(n)$ then $HVE^{R}_{mn}$ is a consistent estimator for $H(f)$.}}
\end{theorem}
{\bf\em {Proof:}} {\textcolor{blue} {Van Es \cite{31me} showed that $HVE_{mn}$ is a consistent estimator of $H(f)$. So
\begin{eqnarray}
\label{1a} E\left[HVE_{mn}\right] &\rightarrow&  H(f),\\
\label{4a} Var\left[HVE_{mn}\right] &\rightarrow& 0,
\end{eqnarray}
as $m=o(n)$.
According to the previous Lemma, $HVE^{R}_{mn}-HVE_{mn}\rightarrow 0$, so (Billingsley \cite{bil})
\begin{eqnarray}
\label{2a}    E\left[HVE^{R}_{mn}-HVE_{mn}\right] &\rightarrow&  0,\\
\label{3a}    E\left[HVE^{R}_{mn}-HVE_{mn}\right]^2 &\rightarrow&  0.
\end{eqnarray}
Now, by using (\ref{1a}) and (\ref{2a}), we conclude $E[HVE^{R}_{mn}]\rightarrow H(f)$.\\
On the other hand, using (\ref{2a}) and (\ref{3a}), we have
\begin{equation}
\label{5a}
Var[HVE^{R}_{mn}-HVE_{mn}]=E[HVE^{R}_{mn}-HVE_{mn}]^2-E^2[HVE^{R}_{mn}-HVE_{mn}]\rightarrow0.
\end{equation}
By applying (\ref{2a}) and (\ref{5a}), we deduce $HVE^{R}_{mn}-HVE_{mn}=o_P(1)$. Because $HVE_{mn}=H(f)+o_P(1)$, so ${HVE^{R}_{mn}=H(f)+o_P(1)}.~\blacksquare$}}


\subsection{RMSE and efficiency comparisons of entropy estimators}
\par
In this subsection, we report the results of a simulation study which compares the performances of the introduced entropy estimators with the estimators proposed by Vasicek \cite{32me}, Van Es \cite{31me}, Ebrahimi {\em et al.} \cite{12me}, Correa \cite{8me}, Wieczorkowski and Grzegorewski \cite{33me} and Zamanzadeh and Arghami \cite{36me} in terms of their SDs and RMSEs. For selected values of $n$, $N = 10000$ samples of size $n$ were generated from Exponential, Normal and Uniform distributions which are the same three distributions considered by Ebrahimi {\em et al.} \cite{12me} and Correa \cite{8me}. Still an open problem in entropy estimation is the optimal choice of $m$ for given $n$. We choose to use the following heuristic formula \cite{17me} as $m=[\sqrt{n}+0.5].$
\par
Tables 1-3 contain the RMSE and SD values of the nine estimators at different sample sizes for each of the three considered distributions. We observe that the proposed estimators perform well compared with other estimators under different distributions. Also, the first estimator, $HVE^{R}_{mn}$ has the least RMSE among the entropy estimators under Exponential and Normal distributions. In addition, under Uniform distribution it behaves better than Vasicek's  estimator, Van Es's estimator, Correa's estimator and Ebrahimi's estimator, four of the most important estimators of entropy. Moreover, the second estimator, $HW^{R}_{mn}$ has the least RMSE among the entropy estimators under Exponential, Normal and Uniform distributions. In many cases, we see that the proposed estimators have the lowest SD among all entropy estimators. Moreover, in all cases, RMSE and SD decrease with the sample size; and $HV_{mn}$ and $HW_{mn}$ have the same SD, because Wieczorkowski and Grzegorzewski modified the Vasicek estimator by adding a bias correction.

\begin{sidewaystable}
\label{bia1}
\caption{{\footnotesize Root mean squared error and standard deviation of the estimators of the entropy of the Exponential distribution.}}
\vspace{0.1in}
\centering
\begin{tabular}{ccccccccccc}
\hline
\hline
\multicolumn{11}{c}{RMSE(SD)} \\
\cline{3-11}
\scriptsize $n$ &\scriptsize $m$ &\scriptsize ${HV}_{mn}$ &\scriptsize ${HVE}_{mn}$ &\scriptsize ${HC}_{mn}$ &\scriptsize ${HE}_{mn}$ &\scriptsize ${HW}_{mn}$ &\scriptsize ${HZ1}_{mn}$ &\scriptsize ${HZ2}_{mn}$ &\scriptsize ${HVE}_{mn}^{R}$ &\scriptsize ${HW}^R_{mn}$ \\
\hline
\scriptsize5  &\scriptsize 2 &\scriptsize 0.930(0.559) &\scriptsize 0.596(0.586) &\scriptsize 0.743(0.554) &\scriptsize 0.651(0.559) &\scriptsize 0.561(0.559) & \scriptsize 0.575(0.573) & \scriptsize 0.576(0.575) &\scriptsize 0.428(0.413)  & \scriptsize 0.425(0.399)\\
\scriptsize10 &\scriptsize 3 &\scriptsize 0.570(0.360) &\scriptsize 0.392(0.373) &\scriptsize 0.435(0.361) &\scriptsize 0.404(0.360) &\scriptsize 0.361(0.360) & \scriptsize 0.391(0.383) & \scriptsize 0.389(0.383) &\scriptsize 0.204(0.202)  & \scriptsize 0.247(0.218)\\
\scriptsize15 &\scriptsize 4 &\scriptsize 0.421(0.284) &\scriptsize 0.310(0.290) &\scriptsize 0.328(0.290) &\scriptsize 0.308(0.284) &\scriptsize 0.282(0.282) & \scriptsize 0.327(0.306) & \scriptsize 0.321(0.305) &\scriptsize 0.133(0.132)  & \scriptsize 0.190(0.158)\\
\scriptsize20 &\scriptsize 4 &\scriptsize 0.356(0.242) &\scriptsize 0.274(0.250) &\scriptsize 0.272(0.247) &\scriptsize 0.263(0.242) &\scriptsize 0.242(0.242) & \scriptsize 0.300(0.261) & \scriptsize 0.286(0.260) &\scriptsize 0.104(0.102)  & \scriptsize 0.150(0.121)\\
\scriptsize30 &\scriptsize 5 &\scriptsize 0.276(0.198) &\scriptsize 0.227(0.201) &\scriptsize 0.208(0.197) &\scriptsize 0.201(0.187) &\scriptsize 0.198(0.198) & \scriptsize 0.266(0.209) & \scriptsize 0.245(0.208) &\scriptsize 0.078(0.069)  & \scriptsize 0.119(0.088)\\
\scriptsize50 &\scriptsize 7 &\scriptsize 0.198(0.150) &\scriptsize 0.181(0.151) &\scriptsize 0.156(0.151) &\scriptsize 0.153(0.150) &\scriptsize 0.150(0.148) & \scriptsize 0.242(0.160) & \scriptsize 0.210(0.158) &\scriptsize 0.066(0.042)  & \scriptsize 0.093(0.061)\\
\hline
\end{tabular}
\vspace{.1in}
\label{bia2}
\caption{{\footnotesize Root mean squared error and standard deviation of the estimators of the entropy of the standard Normal distribution.}}
\vspace{0.1in}
\centering
\begin{tabular}{ccccccccccc}
\hline
\hline
\multicolumn{11}{c}{RMSE(SD)} \\
\cline{3-11}
\scriptsize $n$ &\scriptsize $m$ &\scriptsize ${HV}_{mn}$ &\scriptsize ${HVE}_{mn}$ &\scriptsize ${HC}_{mn}$ &\scriptsize ${HE}_{mn}$ &\scriptsize ${HW}_{mn}$ &\scriptsize ${HZ1}_{mn}$ &\scriptsize ${HZ2}_{mn}$ &\scriptsize ${HVE}_{mn}^{R}$ &\scriptsize ${HW}^R_{mn}$ \\
\hline
\scriptsize5  &\scriptsize 2 &\scriptsize 0.994(0.425) &\scriptsize 0.509(0.452) &\scriptsize 0.793(0.418) &\scriptsize 0.666(0.425) &\scriptsize 0.464(0.413) & \scriptsize 0.494(0.407) & \scriptsize 0.493(0.407) &\scriptsize 0.384(0.370) & \scriptsize 0.355(0.344)\\
\scriptsize10 &\scriptsize 3 &\scriptsize 0.618(0.269) &\scriptsize 0.366(0.283) &\scriptsize 0.407(0.271) &\scriptsize 0.408(0.269) &\scriptsize 0.297(0.264) & \scriptsize 0.303(0.255) & \scriptsize 0.310(0.255) &\scriptsize 0.215(0.181) & \scriptsize 0.197(0.194)\\
\scriptsize15 &\scriptsize 4 &\scriptsize 0.474(0.211) &\scriptsize 0.318(0.220) &\scriptsize 0.348(0.213) &\scriptsize 0.294(0.211) &\scriptsize 0.233(0.211) & \scriptsize 0.222(0.193) & \scriptsize 0.232(0.192) &\scriptsize 0.182(0.121) & \scriptsize 0.147(0.147)\\
\scriptsize20 &\scriptsize 4 &\scriptsize 0.373(0.179) &\scriptsize 0.276(0.185) &\scriptsize 0.265(0.182) &\scriptsize 0.247(0.179) &\scriptsize 0.190(0.178) & \scriptsize 0.190(0.170) & \scriptsize 0.205(0.169) &\scriptsize 0.153(0.095) & \scriptsize 0.118(0.117)\\
\scriptsize30 &\scriptsize 5 &\scriptsize 0.282(0.144) &\scriptsize 0.243(0.148) &\scriptsize 0.194(0.146) &\scriptsize 0.186(0.144) &\scriptsize 0.150(0.144) & \scriptsize 0.148(0.135) & \scriptsize 0.165(0.135) &\scriptsize 0.140(0.065) & \scriptsize 0.091(0.087)\\
\scriptsize50 &\scriptsize 7 &\scriptsize 0.199(0.110) &\scriptsize 0.212(0.110) &\scriptsize 0.134(0.112) &\scriptsize 0.128(0.110) &\scriptsize 0.110(0.109) & \scriptsize 0.108(0.104) & \scriptsize 0.127(0.104) &\scriptsize 0.137(0.039) & \scriptsize 0.070(0.061)\\
\hline
\end{tabular}
\vspace{.1in}
\label{bia3}
\caption{{\footnotesize Root mean squared error and standard deviation of the estimators of the entropy of the Uniform distribution.}}
\vspace{0.1in}
\centering
\begin{tabular}{ccccccccccc}
\hline
\hline
\multicolumn{11}{c}{RMSE(SD)} \\
\cline{3-11}
\scriptsize $n$ &\scriptsize $m$ &\scriptsize ${HV}_{mn}$ &\scriptsize ${HVE}_{mn}$ &\scriptsize ${HC}_{mn}$ &\scriptsize ${HE}_{mn}$ &\scriptsize ${HW}_{mn}$ &\scriptsize ${HZ1}_{mn}$ &\scriptsize ${HZ2}_{mn}$ &\scriptsize ${HVE}_{mn}^{R}$ &\scriptsize ${HW}^R_{mn}$ \\
\hline
\scriptsize5  &\scriptsize 2 &\scriptsize 0.774(0.346) &\scriptsize 0.407(0.407) &\scriptsize 0.566(0.336) &\scriptsize 0.405(0.446) & \scriptsize 0.346(0.346) &\scriptsize 0.330(0.326) &\scriptsize 0.330(0.327) &\scriptsize 0.328(0.297)  & \scriptsize 0.283(0.255)\\
\scriptsize10 &\scriptsize 3 &\scriptsize 0.455(0.167) &\scriptsize 0.216(0.216) &\scriptsize 0.295(0.169) &\scriptsize 0.235(0.167) & \scriptsize 0.166(0.166) &\scriptsize 0.179(0.176)  &\scriptsize 0.180(0.178) &\scriptsize 0.165(0.127)  & \scriptsize 0.134(0.104)\\
\scriptsize15 &\scriptsize 4 &\scriptsize 0.343(0.110) &\scriptsize 0.155(0.155) &\scriptsize 0.208(0.112) &\scriptsize 0.159(0.110) &\scriptsize 0.110(0.110) &\scriptsize 0.137(0.123) & \scriptsize 0.136(0.127) &\scriptsize 0.116(0.082)  & \scriptsize 0.090(0.063)\\
\scriptsize20 &\scriptsize 4 &\scriptsize 0.274(0.087) &\scriptsize 0.126(0.121) &\scriptsize 0.157(0.088) &\scriptsize 0.133(0.087) &\scriptsize 0.087(0.087)  &\scriptsize 0.125(0.100) & \scriptsize 0.121(0.105) &\scriptsize 0.099(0.063)  & \scriptsize 0.078(0.049)\\
\scriptsize30 &\scriptsize 5 &\scriptsize 0.210(0.059) &\scriptsize 0.086(0.086) &\scriptsize 0.110(0.061) &\scriptsize 0.096(0.059) &\scriptsize 0.059(0.059) & \scriptsize 0.112(0.073) & \scriptsize 0.104(0.078) &\scriptsize 0.074(0.041)  & \scriptsize 0.056(0.031)\\
\scriptsize50 &\scriptsize 7 &\scriptsize 0.155(0.037) &\scriptsize 0.057(0.057) &\scriptsize 0.075(0.038) &\scriptsize 0.062(0.037) &\scriptsize 0.037(0.037) & \scriptsize 0.101(0.051) & \scriptsize 0.090(0.055) &\scriptsize 0.052(0.024)  & \scriptsize 0.038(0.017)\\
\hline
\end{tabular}
\end{sidewaystable}

\par
{\textcolor{blue} {Table 4 can be applied for choosing the required number of samples $n$, by imposing some proper fixed values of $m$ and $w$, and preserving  acceptable values of MSE. For example, if a user decides to utilize $HVE_{mn}^{R}$ and required a maximum of $\mbox{MSE}=0.0109$ then, based on the Table 4 assessing  $n\geq 30$, with $m=5$ and $w=5$ would be appropriate. Also, from this table the following general trends can be considered:
\begin{itemize}
{\item In $HVE_{mn}^{R}$ for fixed $w$ and $n$ as $m$ decreases, the MSEs decrease,}
{\item In $HW_{mn}^{R}$ for fixed $w$ and large (small) $n$ as $m$ decreases (increases), the MSEs decrease,}
{\item In $HVE_{mn}^{R}$ for fixed $m$ and $n$ as $w$ increases, the MSEs decrease,}
{\item In $HW_{mn}^{R}$ for fixed $m$ and large (small) $n$ as $w$ decreases (increases), the MSEs decrease.}
\end{itemize}}}
\par
{\textcolor{blue} {We also focus on the efficiency of the entropy estimators based on RSS relative to the entropy estimators based on SRS. However, it is evidently
expected that the entropy estimator based RSS is more efficient than SRS. As the estimators are somehow complicated, we consider the behavior of them by applying the Monte Carlo method and using the idea of Rezakhah and Shemehsavar \cite{Dr}. The values of MSE for estimators are numerically evaluated for different $n$ and some distributions such as Exp(1) and Uniform, see Table 5. Based on the observed values, our conjecture is that for all distributions  $MSE(entropy~ estimator)$ is of order $O(n^{-b})$  and for any fixed distribution there is some constant $C$ that  $MSE(entropy~ estimator)\sim C/n^b$.  We find that for competitor estimators $b=1$ and for proposed estimators $b=1.4$. So our entropy estimators based on RSS is more efficient than the others.}}

\begin{table}
\label{bia1}
\caption{{\footnotesize  Standard deviation and mean squared error of the proposed entropy estimators of the Normal distribution for different $n$, $m$ and $w$, to obtain the sample size.}}
\vspace{0.1in}
\centering
\begin{tabular}{cc|cc|cc||cc|cc}
\hline
\hline
 & \scriptsize$w$ & \multicolumn{4}{c||}{\scriptsize$3$} & \multicolumn{4}{c}{\scriptsize$5$}   \\
 \hline
 &\scriptsize Stat. & \multicolumn{2}{c|}{\scriptsize$HVE_{mn}^{R}$} & \multicolumn{2}{c||}{\scriptsize$HW_{mn}^{R}$} & \multicolumn{2}{c|}{\scriptsize$HVE_{mn}^{R}$} & \multicolumn{2}{c}{\scriptsize$HW_{mn}^{R}$} \\
 \hline
\scriptsize $n$ & \scriptsize$m$ & \scriptsize SD & \scriptsize MSE & \scriptsize SD & \scriptsize MSE & \scriptsize SD & \scriptsize MSE & \scriptsize SD & \scriptsize MSE\\
\hline
 & \scriptsize2 & \scriptsize0.2040 & \scriptsize0.0454 &  \scriptsize0.2106 & \scriptsize0.1240 & \scriptsize 0.2012& \scriptsize 0.0336&  \scriptsize 0.2045& \scriptsize0.0451\\
\scriptsize10 & \scriptsize3  & \scriptsize0.1793 & \scriptsize0.0458  & \scriptsize0.1937 & \scriptsize0.0385 & \scriptsize 0.1741& \scriptsize 0.0395&  \scriptsize 0.1861& \scriptsize0.0350\\
 & \scriptsize4  & \scriptsize0.1705 & \scriptsize0.0467 & \scriptsize0.1862 & \scriptsize0.0383 & \scriptsize 0.1655& \scriptsize 0.0405&  \scriptsize 0.1764& \scriptsize0.0311\\
 & \scriptsize5  & \scriptsize0.1688 & \scriptsize0.0563 & \scriptsize0.1819 & \scriptsize0.0379 & \scriptsize 0.1639& \scriptsize 0.0529&  \scriptsize 0.1728& \scriptsize0.0303\\
\hline
 & \scriptsize3  & \scriptsize0.1277 & \scriptsize0.0267 &  \scriptsize0.1452 & \scriptsize0.0225 & \scriptsize 0.1252& \scriptsize 0.0182&  \scriptsize 0.1407& \scriptsize0.0211\\
\scriptsize15 & \scriptsize4  & \scriptsize0.1211 & \scriptsize0.0333 & \scriptsize0.1426 & \scriptsize0.0221& \scriptsize 0.1166& \scriptsize 0.0194&  \scriptsize 0.1381& \scriptsize0.0200\\
 & \scriptsize5  & \scriptsize0.1192 & \scriptsize0.0384 &  \scriptsize0.1446 & \scriptsize0.0214 & \scriptsize 0.1114& \scriptsize 0.0291&  \scriptsize 0.1397& \scriptsize0.0198\\
 & \scriptsize6  & \scriptsize0.1152 & \scriptsize0.0437 &  \scriptsize0.1442 & \scriptsize0.0204 & \scriptsize 0.1122& \scriptsize 0.0324&  \scriptsize 0.1382& \scriptsize0.0191\\
\hline
 & \scriptsize3  & \scriptsize0.1030 & \scriptsize0.0184 & \scriptsize0.1158 & \scriptsize0.0159 & \scriptsize 0.1003& \scriptsize 0.0104&  \scriptsize 0.1114& \scriptsize0.0147\\
\scriptsize20 & \scriptsize4  & \scriptsize0.0947 & \scriptsize0.0231  & \scriptsize0.1164 & \scriptsize0.0148& \scriptsize 0.0915& \scriptsize 0.0126&  \scriptsize 0.1128& \scriptsize0.0140\\
 & \scriptsize5  &\scriptsize 0.0948 &\scriptsize 0.0235&  \scriptsize0.1184 & \scriptsize0.0148 & \scriptsize 0.0878& \scriptsize 0.0138&  \scriptsize 0.1161& \scriptsize0.0138\\
 & \scriptsize6  & \scriptsize0.0881 &\scriptsize 0.0354 &  \scriptsize0.1210 & \scriptsize0.0147 & \scriptsize 0.0857& \scriptsize 0.0184&  \scriptsize 0.1172& \scriptsize0.0137\\
\hline
 & \scriptsize4 & \scriptsize0.0697 &\scriptsize 0.0148 & \scriptsize0.0839 & \scriptsize0.0079 & \scriptsize 0.0683& \scriptsize 0.0101&  \scriptsize 0.0820& \scriptsize0.0096\\
\scriptsize30 & \scriptsize5 & \scriptsize0.0646 &\scriptsize 0.0197 & \scriptsize0.0858 &\scriptsize 0.0080 & \scriptsize 0.0634& \scriptsize 0.0109&  \scriptsize 0.0833& \scriptsize0.0098\\
 & \scriptsize6 &  \scriptsize0.0626 & \scriptsize0.0281 & \scriptsize0.0886 &\scriptsize 0.0084 & \scriptsize 0.0609& \scriptsize 0.0135&  \scriptsize 0.0871& \scriptsize0.0099\\
 & \scriptsize7 & \scriptsize 0.0603 & \scriptsize0.0310 &  \scriptsize0.0910 & \scriptsize0.0087 & \scriptsize 0.0589& \scriptsize 0.0191&  \scriptsize 0.0895& \scriptsize0.0102\\
\hline
 &\scriptsize 5  & \scriptsize0.0521 & \scriptsize0.0144 & \scriptsize 0.0677 & \scriptsize0.0057 & \scriptsize 0.0508& \scriptsize 0.0096&  \scriptsize 0.0667& \scriptsize0.0073\\
\scriptsize40 &\scriptsize 6  & \scriptsize0.0494 & \scriptsize0.0187  &\scriptsize 0.0707 & \scriptsize0.0060 & \scriptsize 0.0484& \scriptsize 0.0118&  \scriptsize 0.0691& \scriptsize0.0074\\
 &\scriptsize 7  & \scriptsize0.0480 &\scriptsize 0.0190 &  \scriptsize0.0728 & \scriptsize0.0063 & \scriptsize 0.0466& \scriptsize 0.0156&  \scriptsize 0.0708& \scriptsize0.0077\\
 & \scriptsize8  & \scriptsize0.0466 &\scriptsize 0.0283 & \scriptsize 0.0754 & \scriptsize0.0066 & \scriptsize 0.0447& \scriptsize 0.0188&  \scriptsize 0.0722& \scriptsize0.0079\\
\hline
 &\scriptsize 6  &\scriptsize 0.0421 &\scriptsize 0.0148 &\scriptsize 0.0583 & \scriptsize0.0047 & \scriptsize 0.0409& \scriptsize 0.0109&  \scriptsize 0.0572& \scriptsize0.0061\\
\scriptsize50 &\scriptsize 7 &\scriptsize 0.0398 & \scriptsize0.0184 &  \scriptsize0.0607 & \scriptsize0.0049 & \scriptsize 0.0390& \scriptsize 0.0146&  \scriptsize 0.0590& \scriptsize0.0062\\
 & \scriptsize8  &\scriptsize 0.0388 & \scriptsize0.0225 &\scriptsize 0.0634 & \scriptsize0.0054 & \scriptsize 0.0383& \scriptsize 0.0151&  \scriptsize 0.0611& \scriptsize0.0066\\
 & \scriptsize9  & \scriptsize0.0379 & \scriptsize0.0270 &  \scriptsize0.0649 &\scriptsize 0.0057 & \scriptsize 0.0375& \scriptsize 0.0190&  \scriptsize 0.0630& \scriptsize0.0071\\
\hline
 & \scriptsize7  & \scriptsize0.0347 &\scriptsize 0.0154 &  \scriptsize0.0518 & \scriptsize0.0041 & \scriptsize 0.0341& \scriptsize 0.0131&  \scriptsize 0.0509& \scriptsize0.0054\\
\scriptsize60 & \scriptsize8  & \scriptsize0.0334 &\scriptsize 0.0191  &\scriptsize 0.0543 & \scriptsize0.0044 & \scriptsize 0.0330& \scriptsize 0.0142&  \scriptsize 0.0530& \scriptsize0.0057\\
 &\scriptsize 9  & \scriptsize0.0324 &\scriptsize 0.0223 & \scriptsize 0.0569 & \scriptsize0.0049 & \scriptsize 0.0322& \scriptsize 0.0154&  \scriptsize 0.0542& \scriptsize0.0060\\
 \hline
\end{tabular}
\end{table}

\begin{table}
\caption{{\footnotesize Trend of the entropy estimators MSE for different $n$, by considering (a): Exp(1) distribution, (b): Uniform distribution.}}
\centering
\vspace{0.1in}
\begin{tabular}{llllllll}
\hline
\hline
 & & \multicolumn{6}{c}{\scriptsize$n$} \\
\cline{3-8}
 &\scriptsize Dist. & \scriptsize10 & \scriptsize20 & \scriptsize30 & \scriptsize40 & \scriptsize50 & \scriptsize100 \\
\hline
\scriptsize $MSE(HV_{mn})\times n$ & \scriptsize(a)  & \scriptsize2.9435 & \scriptsize2.5908 & \scriptsize 2.1738& \scriptsize2.1382 & \scriptsize1.8911 & \scriptsize 1.8049\\
      & \scriptsize(b)  & \scriptsize 2.0022 & \scriptsize 1.5130& \scriptsize 1.3236& \scriptsize 1.2375& \scriptsize 1.2051& \scriptsize1.2050 \\
\hline
\scriptsize $MSE(HVE_{mn})\times n$ & \scriptsize(a)  & \scriptsize 1.4930& \scriptsize 1.4946& \scriptsize 1.5285& \scriptsize 1.5854& \scriptsize 1.6548& \scriptsize 1.6367\\
      & \scriptsize(b)  & \scriptsize 0.4771 & \scriptsize 0.2995& \scriptsize 0.2180& \scriptsize 0.1891& \scriptsize 0.1765& \scriptsize 0.1706\\
\hline
\scriptsize $MSE(HC_{mn})\times n$ & \scriptsize(a)  & \scriptsize 1.8923& \scriptsize 1.4797& \scriptsize1.2979 & \scriptsize 1.2323& \scriptsize 1.2168& \scriptsize 1.2072\\
      & \scriptsize(b)  & \scriptsize 0.8702 & \scriptsize 0.4930& \scriptsize 0.3630& \scriptsize 0.3043& \scriptsize 0.2813& \scriptsize 0.2800\\
\hline
\scriptsize $MSE(HW_{mn})\times n$ & \scriptsize(a)  & \scriptsize 1.2756& \scriptsize 1.1720& \scriptsize 1.1463& \scriptsize 1.1321& \scriptsize 1.1446& \scriptsize 1.1408\\
      & \scriptsize(b)  & \scriptsize 0.2707 & \scriptsize 0.1487& \scriptsize 0.1056& \scriptsize 0.0842& \scriptsize 0.0676& \scriptsize0.0658 \\
\hline
\scriptsize $MSE(HE_{mn})\times n$ & \scriptsize(a)  & \scriptsize1.5283 & \scriptsize1.4722 & \scriptsize 1.3000& \scriptsize 1.1988& \scriptsize 1.1093& \scriptsize 1.1018\\
      & \scriptsize(b)  & \scriptsize 0.5339 & \scriptsize 0.3658& \scriptsize 0.2808& \scriptsize 0.2442& \scriptsize 0.2096& \scriptsize 0.1900\\
\hline
\scriptsize $MSE(HZ1_{mn})\times n$ & \scriptsize(a)  & \scriptsize 1.5288& \scriptsize 1.8000& \scriptsize 2.1227& \scriptsize 2.5560& \scriptsize 2.9282& \scriptsize 2.9005\\
      & \scriptsize(b)  & \scriptsize 0.3204 & \scriptsize 0.3125& \scriptsize 0.3763& \scriptsize 0.4642& \scriptsize 0.5101& \scriptsize 0.5762\\
\hline
\scriptsize $MSE(HZ2_{mn})\times n$ & \scriptsize(a)  & \scriptsize 1.5132& \scriptsize 1.6359& \scriptsize 1.8007& \scriptsize 2.0880& \scriptsize 2.2050& \scriptsize 2.22302\\
      & \scriptsize(b)  & \scriptsize  0.3240& \scriptsize 0.2928& \scriptsize 0.3245& \scriptsize 0.3750& \scriptsize 0.4050& \scriptsize 0.4333\\
\hline
\scriptsize $MSE(HVE^{R}_{mn})\times n^{1.4}$ & \scriptsize(a)  & \scriptsize 1.0040& \scriptsize 0.7093& \scriptsize 0.7115& \scriptsize 0.7973& \scriptsize 1.0300& \scriptsize 1.0714\\
      & \scriptsize(b)  & \scriptsize 0.7134 & \scriptsize 0.7064& \scriptsize 0.6443& \scriptsize 0.6288& \scriptsize 0.6218& \scriptsize 0.6206\\
\hline
\scriptsize $MSE(HW^{R}_{mn})\times n^{1.4}$ & \scriptsize(a)  & \scriptsize 1.4964& \scriptsize 1.3990& \scriptsize 1.6923& \scriptsize 1.8967& \scriptsize 2.1588& \scriptsize2.1959 \\
      & \scriptsize(b)  & \scriptsize 0.4507 & \scriptsize 0.3824& \scriptsize 0.3759& \scriptsize 0.3679& \scriptsize 0.3599& \scriptsize 0.3579\\
\hline
\end{tabular}
\end{table}

\section{Testing Normality Based on Entropy}
\par
The Normal distribution is the most important distribution in statistics and has a predominant presence in statistical inference. Many statistical techniques are based on the assumption that the data come from this well-known, Bell-shaped distribution. Consequently, the results of these techniques can be completely unreliable if the normality assumption is violated. Thus, it becomes very important to check this assumption in an appropriate and efficient way and that is why goodness-of-fit techniques, especially for Normal distribution, have attracted the attention of many researchers in statistical inference.
\par
In this section, we first introduce a goodness-of-fit test based on the proposed entropy estimators for testing normality and then we compare the powers of the introduced tests with the tests based on the entropy estimators.

\subsection{Test statistics}
\par
A well-known theorem of information theory \cite{30me} states that among all continuous distributions that possess a density function $f$, and have a given variance $\sigma^2$, the entropy $H(f)$ maximized by the Normal distribution. Also simulation results, which provided in some papers such as Wieczorkowski and Grzegorzewsky \cite{33me} and Alizadeh Noughabi \cite{3me}, show that in most cases the test statistic based on $HVE_{mn}$ and $HW_{mn}$ have most power against the alternatives with supports $(-\infty,\infty)$ and $(0,\infty)$, respectively. Based on these properties, following Vasicek \cite{32me}, we introduce the following statistics for testing normality:
\begin{eqnarray*}
TVE^{R}_{mn}&=&\frac{\exp\{HVE^{R}_{mn}\}}{\hat\sigma}=\frac{1}{\hat{\sigma}}\exp\left\{\sum_{k=m}^{n}\frac{1}{k}\right\}\left[\prod_{i=1}^{n-m}
(Y^R_{i+m}-Y^R_{i})\right]^\frac{1}{n-m},\\
TW^{R}_{mn}&=&\frac{\exp\{HW^{R}_{mn}\}}{\hat\sigma}=\frac{1}{\hat{\sigma}}\left[\prod_{i=1}^{n}(Y^R_{i+m}-Y^R_{i-m})\right]^\frac{1}{n}\exp\{c\},\\
\end{eqnarray*}
where $HVE^{R}_{mn}$, $HW^{R}_{mn}$ and $c$ are defined in (\ref{VER}), (\ref{WR}) and (\ref{c}), respectively. \par
Small values of $TVE^{R}_{mn}$ and $TW^{R}_{mn}$ can be regarded as a symptom of non-normality and therefore we reject the hypothesis of normality for small enough values of $TVE^{R}_{mn}$ and $TW^{R}_{mn}$. We see that $TVE^{R}_{mn}$ and $TW^{R}_{mn}$ are invariant with respect to transformations of location and scale.
\par
Because the sampling distributions of the test statistics are complicated, we determine the percentage points using 10000 Monte Carlo samples from a Normal distribution. So, the exact critical values (non-asymptotic) of the test statistics $TVE^{R}_{mn}$ and $TW^{R}_{mn}$, for various sample sizes by Monte Carlo
simulations are given in Tables 6 and 7, respectively.

\begin{table}
\label{CV1}
\caption{{\footnotesize Critical values of the $TVE^{R}_{mn}$ statistic for testing normality with $\alpha=0.05$, for various sample size and different $m$. }}
\vspace{0.2in}
\hspace{-.1in}
\centering
\begin{tabular}{ccccccccccc}
\hline
\hline
\multicolumn{11}{c}{\scriptsize $m$} \\
\cline{2-11}
\scriptsize $n$ &\scriptsize  1 &\scriptsize  2 &\scriptsize 3 &\scriptsize 4 &\scriptsize 5 &\scriptsize 6 &\scriptsize 7 &\scriptsize 8 &\scriptsize 9 &\scriptsize 10 \\
\hline
\scriptsize5   &\scriptsize 2.2978 &\scriptsize 2.8789 &        &        &        &        &        &        &        &  \\
\scriptsize6   &\scriptsize 2.4129 &\scriptsize 2.9539 &\scriptsize 3.0372 &        &        &        &        &        &        &  \\
\scriptsize7   &\scriptsize 2.4425 &\scriptsize 2.9681 &\scriptsize 3.0186 &        &        &        &        &        &        &  \\
\scriptsize8   &\scriptsize 2.6098 &\scriptsize 3.0016 &\scriptsize 3.0298 &\scriptsize 3.0811 &        &        &        &        &        &  \\
\scriptsize9   &\scriptsize 2.6505 &\scriptsize 3.0632 &\scriptsize 3.0717 &\scriptsize 3.0362 &        &        &        &        &        &  \\
\scriptsize10  &\scriptsize 2.7015 &\scriptsize 3.1098 &\scriptsize 3.1204 &\scriptsize 3.0619 &\scriptsize 3.0777 &        &        &        &        &  \\
\scriptsize15  &\scriptsize 2.9185 &\scriptsize 3.2197 &\scriptsize 3.2483 &\scriptsize 3.1809 &\scriptsize 3.1102 &\scriptsize 3.0563 &\scriptsize 3.0521 &        &        &  \\
\scriptsize20  &\scriptsize 3.0864 &\scriptsize 3.3446 &\scriptsize 3.3320 &\scriptsize 3.2758 &\scriptsize 3.2090 &\scriptsize 3.1432 &\scriptsize 3.0925 &\scriptsize 3.0642 &\scriptsize 3.0460 &\scriptsize 3.0518 \\
\scriptsize25  &\scriptsize 3.1749 &\scriptsize 3.4084 &\scriptsize 3.4155 &\scriptsize 3.3585 &\scriptsize 3.2926 &\scriptsize 3.2317 &\scriptsize 3.1718 &\scriptsize 3.1237 &\scriptsize 3.0881 &\scriptsize 3.0584 \\
\scriptsize30  &\scriptsize 3.2583 &\scriptsize 3.4758 &\scriptsize 3.4760 &\scriptsize 3.4202 &\scriptsize 3.3672 &\scriptsize 3.2962 &\scriptsize 3.2427 &\scriptsize 3.1920 &\scriptsize 3.1551 &\scriptsize 3.1093 \\
\scriptsize40  &\scriptsize 3.3616 &\scriptsize 3.5458 &\scriptsize 3.5557 &\scriptsize 3.5069 &\scriptsize 3.4614 &\scriptsize 3.3974 &\scriptsize 3.3495 &\scriptsize 3.2923 &\scriptsize 3.2521 &\scriptsize 3.2092 \\
\scriptsize50  &\scriptsize 3.4557 &\scriptsize 3.6218 &\scriptsize 3.6174 &\scriptsize 3.5827 &\scriptsize 3.5312 &\scriptsize 3.4842 &\scriptsize 3.4275 &\scriptsize 3.3843 &\scriptsize 3.3380 &\scriptsize 3.3018 \\
\hline
\end{tabular}
\end{table}

\begin{table}
\label{CV2}
\caption{{\footnotesize Critical values of the $TW^{R}_{mn}$ statistic for testing normality with $\alpha=0.05$, for various sample size and different $m$.}}
\vspace{0.2in}
\hspace{-.1in}
\centering
\vspace{0.1in}
\begin{tabular}{ccccccccccc}
\hline
\hline
\multicolumn{11}{c}{\scriptsize $m$} \\
\cline{2-11}
\scriptsize $n$ &\scriptsize  1 &\scriptsize  2 &\scriptsize 3 &\scriptsize 4 &\scriptsize 5 &\scriptsize 6 &\scriptsize 7 &\scriptsize 8 &\scriptsize 9 &\scriptsize 10 \\
\hline
\scriptsize5   &\scriptsize 1.2213 &\scriptsize 1.7656 &        &        &        &        &        &        &        &  \\
\scriptsize6   &\scriptsize 1.3840 &\scriptsize 1.8527 &\scriptsize 1.9232 &        &        &        &        &        &        &  \\
\scriptsize7   &\scriptsize 1.5162 &\scriptsize 1.9638 &\scriptsize 2.0650 &        &        &        &        &        &        &  \\
\scriptsize8   &\scriptsize 1.6828 &\scriptsize 2.0872 &\scriptsize 2.1487 &\scriptsize 2.0959 &        &        &        &        &        &  \\
\scriptsize9   &\scriptsize 1.7952 &\scriptsize 2.1781 &\scriptsize 2.2488 &\scriptsize 2.1969 &        &        &        &        &        &  \\
\scriptsize10  &\scriptsize 1.8958 &\scriptsize 2.3111 &\scriptsize 2.3372 &\scriptsize 2.2995 &\scriptsize 2.1925 &        &        &        &        &  \\
\scriptsize15  &\scriptsize 2.2266 &\scriptsize 2.6853 &\scriptsize 2.7315 &\scriptsize 2.6744 &\scriptsize 2.5836 &\scriptsize 2.4902 &\scriptsize 2.3954 &        &        &  \\
\scriptsize20  &\scriptsize 2.4118 &\scriptsize 2.9098 &\scriptsize 2.9832 &\scriptsize 2.9519 &\scriptsize 2.8859 &\scriptsize 2.7923 &\scriptsize 2.7025 &\scriptsize 2.6117 &\scriptsize 2.5274 &\scriptsize 2.4549 \\
\scriptsize25  &\scriptsize 2.5367 &\scriptsize 3.0491 &\scriptsize 3.1491 &\scriptsize 3.1406 &\scriptsize 3.0935 &\scriptsize 3.0117 &\scriptsize 2.9412 &\scriptsize 2.8675 &\scriptsize 2.7838 &\scriptsize 2.7070 \\
\scriptsize30  &\scriptsize 2.6214 &\scriptsize 3.1331 &\scriptsize 3.2684 &\scriptsize 3.2679 &\scriptsize 3.2360 &\scriptsize 3.1814 &\scriptsize 3.1181 &\scriptsize 3.0506 &\scriptsize 2.9737 &\scriptsize 2.9089 \\
\scriptsize40  &\scriptsize 2.7347 &\scriptsize 3.2665 &\scriptsize 3.4081 &\scriptsize 3.4433 &\scriptsize 3.4263 &\scriptsize 3.3975 &\scriptsize 3.3503 &\scriptsize 3.3037 &\scriptsize 3.2598 &\scriptsize 3.2015 \\
\scriptsize50  &\scriptsize 2.8053 &\scriptsize 3.3353 &\scriptsize 3.4917 &\scriptsize 3.5414 &\scriptsize 3.5423 &\scriptsize 3.5207 &\scriptsize 3.4987 &\scriptsize 3.4557 &\scriptsize 3.4176 &\scriptsize 3.3824 \\
\hline
\end{tabular}
\end{table}

\subsection{Power results}
\par
There are some of test statistics for normality concerning uncensored data. We consider the power of the proposed test statistics in three cases.
\begin{description}
  \item[]{\bf Case I: Test statistics based on the entropy estimators} \par
  Some of the test statistics provided testing normality based on entropy. We consider here the test statistics of Vasicek \cite{32me}, Van Es \cite{31me}, Correa \cite{8me}, Wieczorkowski and Grzegorewski \cite{33me} and Zamanzadeh and Arghami \cite{36me} among them. The test statistics proposed by them are
\begin{eqnarray*}
TV_{mn}&=&\frac{\exp\{HV_{mn}\}}{\hat\sigma}=\frac{n}{2m\hat{\sigma}}\left[\prod_{i=1}^{n}(X_{(i+m)}-X_{(i-m)})\right]^{\frac{1}{n}},\\
TVE_{mn}&=&\frac{\exp\{HVE_{mn}\}}{\hat\sigma}=\frac{1}{\hat{\sigma}}\exp\left\{\sum_{k=m}^{n}\frac{1}{k}\right\}\left[\prod_{i=1}^{n-m}(X_{(i+m)}-X_{(i)})\right]
^\frac{1}{n-m},\\
TC_{mn}&=&\frac{\exp\{HC_{mn}\}}{\hat\sigma}=\frac{1}{\hat{\sigma}}\left[\prod_{i=1}^{n}\frac{\sum_{j=i-m}^{i+m}(X_{(j)}-\bar{X}_{(i)})(j-i)}
{n\sum_{j=i-m}^{i+m}(X_{(j)}-\bar{X}_{(i)})^2}\right]^{-\frac{1}{n}},\\
TW_{mn}&=&\frac{\exp\{HW_{mn}\}}{\hat\sigma}=\frac{1}{\hat{\sigma}}\left[\prod_{i=1}^{n}(X_{(i+m)}-X_{(i-m)})\right]^\frac{1}{n}\exp\{c\},\\
TZ1_{mn}&=&\frac{\exp\{HZ1_{mn}\}}{\hat\sigma}=\frac{1}{\hat{\sigma}}\left[\prod_{i=1}^{n}b_i\right]^\frac{1}{n},
TZ2_{mn}=\frac{\exp\{HZ2_{mn}\}}{\hat\sigma}=\frac{1}{\hat{\sigma}}\left[\prod_{i=1}^{n}b_i\right]^{w_i},
\end{eqnarray*}
where $\bar{X}_{(i)}$, $c$, $b_i$ and $w_i$ are defined in (\ref{xbar}), (\ref{c}), (\ref{bi}) and (\ref{wi}), respectively.
\par
For power comparisons, we compute the powers of the tests based on statistics $TV_{mn}$, $TVE_{mn}$, $TC_{mn}$, $TW_{mn}$, $TZ1_{mn}$, $TZ2_{mn}$, $TVE^{R}_{mn}$ and $TW^{R}_{mn}$ by means of Monte Carlo simulations under $24$ alternatives. Since Wieczorkowski and Grzegorzewski modified the Vasicek estimator by adding a bias correction and we observed that the estimators $HV_{mn}$ and $HW_{mn}$ have the same standard deviation, the tests based on $TV_{mn}$ and $TW_{mn}$ statistics have the same power. The alternatives can be classified into four groups, according to their supports and the shapes of their densities. From the point of view of applied statistics, natural alternatives to Normal distribution are in groups I and II. For the sake of completeness, we also consider groups III and IV. This gives additional insight in understanding the behavior of the new test statistics $TVE^{R}_{mn}$ and $TW^{R}_{mn}$. Esteban {\em et al.} \cite{14me}, in their study of power comparisons of several tests for normality, suggested classifying the alternatives into the following four groups:
\par
\begin{itemize}
\item[] {\em Group I: Support ($-\infty,\infty$), symmetric,} which are presented in Table 9.

\par
\item[]{\em Group II: Support ($-\infty,\infty$), asymmetric,} which are presented in Table 10.

\par
\item[]{\em Group III: Support ($0,1$),} which are presented in Table 11.

\par
\item[]{\em Group IV: Support ($0,\infty$),} which are presented in Table 12.

\end{itemize}

Unfortunately, the powers of the proposed tests depend on the windows size $m$ and the alternative distribution and therefore it is not possible to determine the best value of $m$, for which the tests attain maximum power under all alternatives. Therefore, we use the values of $m$ for which the aforementioned tests attain good (not best) powers for all alternative distributions. These values of $m$ are tabulated in Table 8. The simulation results in \cite{36me} showed that $TZ1_{mn}$ and $TZ2_{mn}$ have not any primacy to others against the alternatives in group II, III and IV. Therefore, we only use them in group I.

\begin{table}
\label{C.V}
\caption{{\footnotesize Proposed windows size $m$ for testing normality for different values of $n$.}}
\centering
\vspace{0.1in}
\begin{tabular}{llllll}
\hline
\hline
\scriptsize $n$ &&&&\scriptsize $m$ \\
\hline
\scriptsize$n\leq8$          &&&&\scriptsize  1\\
\scriptsize$9\leq n\leq 15$   &&&&\scriptsize  2\\
\scriptsize$16\leq n\leq 35$  &&&&\scriptsize  3\\
\scriptsize$36\leq n\leq 60$  &&&&\scriptsize  4\\
\scriptsize$61\leq n\leq 80$  &&&&\scriptsize  5\\
\scriptsize$81\leq n\leq 100$ &&&&\scriptsize  6\\
\hline
\end{tabular}
\end{table}

\par
Tables 9-12 contain the results of $10000$ simulations (of samples size $10$, $20$, $40$ and $50$) per case to obtain the power of the proposed tests and those of the competitor tests, at a significance level $\alpha=0.05$. For each sample size and alternative, the bold type in these tables indicates the statistics achieving the maximal power.

\begin{sidewaystable}
\label{g3}
\caption{{\scriptsize Power comparisons based on different statistics and several sample sizes under the alternatives in group I with $\alpha=0.05$.}}
\vspace{0.1in}
\centering
\begin{tabular}{c|l|cccccc||cccccc|c}
\hline
\hline
\scriptsize $n$ &\scriptsize Alt &\scriptsize $TV_{mn}$ &\scriptsize $TVE_{mn}$ &\scriptsize $TC_{mn}$  &\scriptsize $TZ1_{mn}$  &\scriptsize $TZ2_{mn}$   &\scriptsize $TVE^{R}_{mn}$ &\scriptsize  $TV_{mn}$ &\scriptsize $TVE_{mn}$ &\scriptsize $TC_{mn}$  &\scriptsize $TZ1_{mn}$  &\scriptsize $TZ2_{mn}$  &\scriptsize $TVE^{R}_{mn}$ &\scriptsize  n\\
\hline
     &\scriptsize t(1)           &\scriptsize 0.4360 &\scriptsize 0.6084 &\scriptsize 0.3989 &\scriptsize 0.6320 &\scriptsize 0.6380 & \scriptsize \bf 0.6958&\scriptsize 0.9607 &\scriptsize 0.9877 &\scriptsize 0.9488  &\scriptsize 0.9880 &\scriptsize 0.9930 & \scriptsize \bf 1.0000&\\
     \scriptsize10   &\scriptsize t(3)           &\scriptsize 0.1011 &\scriptsize 0.1877 &\scriptsize 0.0860  &\scriptsize 0.2120 &\scriptsize 0.2160 & \scriptsize \bf 0.2190 &\scriptsize 0.2952 &\scriptsize 0.5567 &\scriptsize 0.2492 &\scriptsize0.5610 &\scriptsize 0.6220& \scriptsize \bf 0.7450&\scriptsize 40\\
     &\scriptsize DE             &\scriptsize 0.0693 &\scriptsize 0.1686 &\scriptsize 0.0587 &\scriptsize 0.1770 &\scriptsize 0.1810 &\scriptsize \bf 0.1810
     &\scriptsize 0.1928 &\scriptsize 0.4533 &\scriptsize 0.1627 &\scriptsize0.4760 &\scriptsize 0.5680 &\scriptsize \bf 0.6670&\\
     &\scriptsize Logistic       &\scriptsize 0.0517 &\scriptsize 0.0915 &\scriptsize 0.0511 &\scriptsize 0.0890 &\scriptsize 0.0910 & \scriptsize \bf 0.0974
     &\scriptsize 0.0550 &\scriptsize 0.1715 &\scriptsize 0.0488 &\scriptsize0.1660 &\scriptsize 0.2110&\scriptsize \bf 0.2280&\\
\hline
     &\scriptsize t(1)           &\scriptsize 0.7296 &\scriptsize 0.8719 &\scriptsize 0.6984 &\scriptsize 0.8850 &\scriptsize 0.9000 & \scriptsize \bf 0.9687& \scriptsize 0.9884 &\scriptsize 0.9968 &\scriptsize 0.9832 &\scriptsize0.9910 &\scriptsize  0.9960&\scriptsize \bf 1.0000& \\
 \scriptsize20  &\scriptsize t(3)          &\scriptsize 0.1571 &\scriptsize 0.3227 &\scriptsize 0.1282  &\scriptsize 0.3770 &\scriptsize 0.4020 &\scriptsize \bf 0.4597 &\scriptsize 0.3629 &\scriptsize 0.6399 &\scriptsize 0.3346  &\scriptsize0.6150 &\scriptsize0.6950 &\scriptsize \bf 0.8148&\scriptsize 50\\
     &\scriptsize DE             &\scriptsize 0.0877 &\scriptsize 0.2655 &\scriptsize 0.0755 &\scriptsize 0.3090 &\scriptsize 0.3440 &\scriptsize \bf 0.3778& \scriptsize 0.2922 &\scriptsize 0.5533 &\scriptsize 0.2419 &\scriptsize0.5540 &\scriptsize 0.6510& \scriptsize \bf 0.7230& \\
     &\scriptsize Logistic       &\scriptsize 0.0454 &\scriptsize 0.1097 &\scriptsize 0.0408 &\scriptsize 0.1330 &\scriptsize 0.1470 &\scriptsize \bf 0.1526& \scriptsize 0.0647 &\scriptsize 0.1976 &\scriptsize 0.0552  &\scriptsize 0.2060 &\scriptsize 0.2480 &\scriptsize \bf 0.2564& \\
\hline
\end{tabular}
\vspace{.1in}

\label{g3}
\caption{{\scriptsize Power comparisons based on different statistics and several sample sizes under the alternatives in group II with $\alpha=0.05$.}}
\vspace{0.1in}
\centering
\begin{tabular}{c|l|cccc||cccc|c}
\hline
\hline
\scriptsize $n$ &\scriptsize Alt &\scriptsize $TV_{mn}$ &\scriptsize $TVE_{mn}$ &\scriptsize $TC_{mn}$   &\scriptsize $TVE^{R}_{mn}$ &\scriptsize  $TV_{mn}$ &\scriptsize $TVE_{mn}$ &\scriptsize $TC_{mn}$   &\scriptsize $TVE^{R}_{mn}$ &\scriptsize  n\\
\hline
     &\scriptsize EV(0,1)           &\scriptsize 0.1038 &\scriptsize 0.1236 &\scriptsize  0.0961  &\scriptsize \bf 0.1240&\scriptsize 0.4008 &\scriptsize 0.3729 &\scriptsize  0.3871  &\scriptsize \bf 0.5250 &\\
\scriptsize10&\scriptsize EV(0,2)           &\scriptsize 0.1036 &\scriptsize 0.1285 &\scriptsize 0.1001   &\scriptsize \bf 0.1383&\scriptsize 0.4023 &\scriptsize 0.3608 &\scriptsize 0.3793   &\scriptsize \bf 0.5241&\scriptsize 40\\
   &\scriptsize EV(0,3)           &\scriptsize 0.1018 &\scriptsize 0.1290 &\scriptsize 0.0944  &\scriptsize \bf 0.1423&\scriptsize 0.4041 &\scriptsize 0.3672 &\scriptsize 0.3876  &\scriptsize \bf 0.5301&\\
     &\scriptsize EV(0,4)           &\scriptsize 0.1075 &\scriptsize 0.1258 &\scriptsize 0.0942  &\scriptsize \bf 0.1410&\scriptsize 0.4051 &\scriptsize 0.3597 &\scriptsize 0.3933   &\scriptsize \bf 0.5227&\\
\hline
     &\scriptsize EV(0,1)           &\scriptsize 0.1886 &\scriptsize 0.1985 &\scriptsize  0.1903  &\scriptsize \bf 0.2779 &\scriptsize 0.4771 &\scriptsize 0.4392 &\scriptsize  0.4853  &\scriptsize \bf 0.6227 &\\
\scriptsize20&\scriptsize EV(0,2)           &\scriptsize 0.1952 &\scriptsize 0.1971 &\scriptsize 0.1933   &\scriptsize \bf 0.2822& \scriptsize 0.4718 &\scriptsize 0.4430 &\scriptsize 0.4871   &\scriptsize \bf 0.6221& \scriptsize 50\\
   &\scriptsize EV(0,3)           &\scriptsize 0.1981 &\scriptsize 0.1969 &\scriptsize 0.1878   &\scriptsize \bf 0.2793&\scriptsize 0.4670 &\scriptsize 0.4344 &\scriptsize 0.4860   &\scriptsize \bf 0.6330&\\
     &\scriptsize EV(0,4)           &\scriptsize 0.1950 &\scriptsize 0.1991 &\scriptsize 0.1937  &\scriptsize \bf 0.2750& \scriptsize 0.4776 &\scriptsize 0.4460 &\scriptsize 0.4792  &\scriptsize \bf 0.6342& \\
\hline
\end{tabular}
\vspace{.1in}

\label{g3}
\caption{{\scriptsize Power comparisons based on different statistics and several sample sizes under the alternatives in group III with $\alpha=0.05$.}}
\vspace{0.1in}
\centering
\begin{tabular}{c|l|cccc||cccc|c}
\hline
\hline
\scriptsize $n$ &\scriptsize Alt &\scriptsize $TV_{mn}$ &\scriptsize $TVE_{mn}$ &\scriptsize $TC_{mn}$  &\scriptsize $TW^{R}_{mn}$ &\scriptsize  $TV_{mn}$ &\scriptsize $TVE_{mn}$ &\scriptsize $TC_{mn}$  &\scriptsize $TW^{R}_{mn}$  &\scriptsize  n\\
\hline
     & \tiny Beta(3,0.5)           &\scriptsize 0.6771 &\scriptsize 0.5560 &\scriptsize  0.6657 &\scriptsize \bf 0.8341 &\scriptsize 0.9998 &\scriptsize 0.9977 &\scriptsize 0.9999 &\scriptsize \bf 1.0000 & \\
   \scriptsize10& \tiny Beta(0.5,3)         &\scriptsize 0.6831 &\scriptsize 0.5647 &\scriptsize 0.6660  &\scriptsize \bf 0.8549 &\scriptsize 1.0000 &\scriptsize 0.9981 &\scriptsize 1.0000  &\scriptsize \bf 1.0000 &\scriptsize 40\\
     & \tiny Beta(0.5,1)         &\scriptsize 0.4784 &\scriptsize 0.2973 &\scriptsize 0.4731  &\scriptsize \bf 0.4967 &\scriptsize 0.9990 &\scriptsize 0.8922 &\scriptsize 0.9984  &\scriptsize \bf 1.0000 & \\
\hline
     & \tiny Beta(3,0.5)           &\scriptsize 0.9796 &\scriptsize 0.8767 &\scriptsize  0.9732 &\scriptsize \bf 1.0000 &\scriptsize 1.0000 &\scriptsize 0.9998 &\scriptsize 1.0000 &\scriptsize \bf 1.0000 & \\
   \scriptsize20& \tiny Beta(0.5,3)         &\scriptsize 0.9761 &\scriptsize 0.8685 &\scriptsize 0.9701  &\scriptsize \bf 1.0000 &\scriptsize 1.0000 &\scriptsize 0.9995 &\scriptsize 1.0000  &\scriptsize \bf 1.0000 &\scriptsize 50\\
     & \tiny Beta(0.5,1)         &\scriptsize 0.8859 &\scriptsize 0.5198 &\scriptsize 0.8702  &\scriptsize \bf 0.9995 &\scriptsize 0.9999 &\scriptsize 0.9465 &\scriptsize {0.9999}  &\scriptsize \bf 1.0000 & \\
\hline
\end{tabular}
\end{sidewaystable}

\begin{sidewaystable}
\label{g3}
\caption{{\scriptsize Power comparisons based on different statistics and several sample sizes under the alternatives in group IV with $\alpha=0.05$.}}
\vspace{0.1in}
\centering
\begin{tabular}{c|l|cccc||cccc|c}
\hline
\hline
\scriptsize $n$ &\scriptsize Alt &\scriptsize $TV_{mn}$ &\scriptsize $TVE_{mn}$ &\scriptsize $TC_{mn}$ &\scriptsize $TW^{R}_{mn}$ &\scriptsize  $TV_{mn}$ &\scriptsize $TVE_{mn}$ &\scriptsize $TC_{mn}$ &\scriptsize $TW^{R}_{mn}$ & \scriptsize n\\
\hline
     &\scriptsize Chi(1)           &\scriptsize 0.7894 &\scriptsize 0.6920 &\scriptsize 0.7748  &\scriptsize \bf 0.9316 &\scriptsize 1.0000 &\scriptsize 0.9999 &\scriptsize 1.0000 &\scriptsize \bf 1.0000 & \\
     &\scriptsize Chi(3)           &\scriptsize 0.2637 &\scriptsize 0.2365 &\scriptsize 0.2527  &\scriptsize \bf 0.3139 &\scriptsize 0.9465 &\scriptsize 0.7792 &\scriptsize 0.9363  &\scriptsize \bf 1.0000 & \\
     &\scriptsize Chi(4)           &\scriptsize 0.1875 &\scriptsize 0.1810 &\scriptsize 0.1839  &\scriptsize \bf 0.2190 & \scriptsize 0.8246 &\scriptsize 0.6209 &\scriptsize 0.8146  &\scriptsize \bf 0.9992 & \\
     &\scriptsize LN(0,0.6)         &\scriptsize 0.2440 &\scriptsize 0.2573 &\scriptsize 0.2350  &\scriptsize \bf 0.3026 &\scriptsize 0.8845 &\scriptsize 0.7821 &\scriptsize 0.8735  &\scriptsize \bf 0.9993 & \\
     &\scriptsize LN(0,1)          &\scriptsize 0.5700 &\scriptsize 0.5159 & \scriptsize0.5544  &\scriptsize \bf 0.7419 &\scriptsize 0.9991 &\scriptsize 0.9903 &\scriptsize 0.9984  &\scriptsize \bf 1.0000 & \\
     &\scriptsize LN(0,1.2)        &\scriptsize 0.7013 &\scriptsize 0.6296 &\scriptsize 0.6908  &\scriptsize \bf 0.8715 &\scriptsize 1.0000 &\scriptsize 0.9984 &\scriptsize 1.0000  &\scriptsize \bf 1.0000 & \\
   \scriptsize10&\scriptsize Gamma(0.5)        &\scriptsize 0.7861 &\scriptsize 0.6937 &\scriptsize 0.7735  &\scriptsize \bf 0.9299 &\scriptsize 1.0000 &\scriptsize 0.9999 &\scriptsize 1.0000  &\scriptsize \bf 1.0000 &\scriptsize 40\\
     &\scriptsize Gamma(1.5)       &\scriptsize 0.2524 &\scriptsize 0.2298 &\scriptsize 0.2546 &\scriptsize \bf 0.3160 &\scriptsize 0.9440 &\scriptsize 0.7902 &\scriptsize 0.9380  &\scriptsize \bf 1.0000 & \\
     &\scriptsize Gamma(2)         &\scriptsize 0.1899 &\scriptsize 0.1831 &\scriptsize 0.1748  &\scriptsize \bf 0.2149 &\scriptsize 0.8273 &\scriptsize 0.6269 &\scriptsize 0.8168  &\scriptsize \bf 0.9996 & \\
     &\scriptsize Weibull(0.5)      &\scriptsize 0.9323 &\scriptsize 0.8857 &\scriptsize 0.9282  &\scriptsize \bf 0.9913 &\scriptsize 1.0000 &\scriptsize 1.0000 &\scriptsize 1.0000  &\scriptsize \bf 1.0000 & \\
     &\scriptsize Weibull(2)       &\scriptsize 0.0800 & \scriptsize0.0774 &\scriptsize 0.0787  &\scriptsize \bf 0.0893 &\scriptsize 0.2633 &\scriptsize 0.1367 &\scriptsize 0.2615 &\scriptsize \bf 0.5042 &  \\
     &\scriptsize Exp(1)           &\scriptsize 0.4249 &\scriptsize 0.3545 &\scriptsize 0.4110  &\scriptsize \bf 0.5416 &\scriptsize 0.9965 &\scriptsize 0.9491 &\scriptsize 0.9949  &\scriptsize \bf 1.0000 & \\
     &\scriptsize Gexp(0.5)         &\scriptsize 0.7765 &\scriptsize 0.6663 &\scriptsize 0.7505  &\scriptsize \bf 0.9290 &\scriptsize 1.0000 &\scriptsize 0.9998 &\scriptsize 1.0000  &\scriptsize \bf 1.0000 & \\

\hline
     &\scriptsize Chi(1)           &\scriptsize 0.9928 &\scriptsize 0.9557 &\scriptsize 0.9895  &\scriptsize \bf 1.0000 &\scriptsize 1.0000 &\scriptsize 1.0000 &\scriptsize 1.0000  &\scriptsize \bf 1.0000 & \\
     &\scriptsize Chi(3)           &\scriptsize 0.6145 &\scriptsize 0.4321 &\scriptsize 0.6067  &\scriptsize \bf 0.9300 &\scriptsize 0.9808 &\scriptsize 0.8749 &\scriptsize 0.9795  &\scriptsize \bf 1.0000 & \\
     &\scriptsize Chi(4)           &\scriptsize 0.4494 &\scriptsize 0.3234 &\scriptsize 0.4438  &\scriptsize \bf 0.7580 &\scriptsize 0.9041 &\scriptsize 0.7251 &\scriptsize 0.9097  &\scriptsize \bf 1.0000 & \\
     &\scriptsize LN(0,0.6)         &\scriptsize 0.5430 &\scriptsize 0.4598 &\scriptsize 0.5403  &\scriptsize \bf 0.8489 &\scriptsize 0.9371 &\scriptsize 0.8707 &\scriptsize 0.9390  &\scriptsize \bf 1.0000 & \\
     &\scriptsize LN(0,1)          &\scriptsize 0.9248 &\scriptsize 0.8340 &\scriptsize 0.9157  &\scriptsize \bf 0.9994 &\scriptsize 1.0000 &\scriptsize 0.9981 &\scriptsize 1.0000  &\scriptsize \bf 1.0000 & \\
     &\scriptsize LN(0,1.2)        &\scriptsize 0.9716 &\scriptsize0.9157 &\scriptsize 0.9695  &\scriptsize \bf 1.0000 &\scriptsize 1.0000 &\scriptsize 0.9999 & \scriptsize0.9999  &\scriptsize \bf 1.0000 & \\
   \scriptsize20&\scriptsize Gamma(0.5)        &\scriptsize 0.9927 &\scriptsize 0.9526 &\scriptsize 0.9895  &\scriptsize \bf 1.0000 &\scriptsize 1.0000 &\scriptsize 1.0000 &\scriptsize 1.0000  &\scriptsize \bf 1.0000 &\scriptsize 50 \\
     &\scriptsize Gamma(1.5)       &\scriptsize 0.6124 &\scriptsize 0.4432 &\scriptsize 0.6052 &\scriptsize \bf 0.9335&\scriptsize 0.9793 &\scriptsize 0.8702 &\scriptsize 0.9778  &\scriptsize \bf 1.0000 & \\
     &\scriptsize Gamma(2)         &\scriptsize 0.4590 &\scriptsize 0.3228 &\scriptsize 0.4464  &\scriptsize \bf 0.7612 &\scriptsize 0.9082 &\scriptsize 0.7214 &\scriptsize 0.9073  &\scriptsize \bf 1.0000& \\
     &\scriptsize Weibull(0.5)      &\scriptsize 0.9994 &\scriptsize 0.9957 &\scriptsize 0.9997  &\scriptsize \bf 1.0000 &\scriptsize 1.0000 &\scriptsize 1.0000 &\scriptsize 1.0000  &\scriptsize \bf 1.0000 & \\
     &\scriptsize Weibull(2)       &\scriptsize 0.1255 &\scriptsize 0.0817 &\scriptsize 0.1256  &\scriptsize \bf 0.1840 &\scriptsize 0.3133 &\scriptsize 0.1529 &\scriptsize 0.3393  &\scriptsize \bf 0.6562 &  \\
     &\scriptsize Exp(1)           &\scriptsize 0.8419 &\scriptsize 0.6557 & \scriptsize0.8330  & \scriptsize\bf 0.9987 & \scriptsize0.9994 &\scriptsize 0.9814 &\scriptsize 0.9994  &\scriptsize \bf 1.0000 & \\
     &\scriptsize Gexp(0.5,1)       & \scriptsize0.9928 &\scriptsize 0.9492 &\scriptsize 0.9891  &\scriptsize \bf 1.0000 &\scriptsize 1.0000 &\scriptsize 1.0000 &\scriptsize 1.0000  &\scriptsize \bf 1.0000 & \\
\hline
\end{tabular}
\end{sidewaystable}
\par
From Tables 9-12, we can see that the tests compared considerably differ in powers. Also, in all cases, the power increases with the sample size.
In groups I and II, it is seen that the tests $TVE^{R}_{mn}$ has the most power and the test $TC_{mn}$ has the least power. Also, the difference of powers between $TV_{mn}$ and $TC_{mn}$ is small. But the difference of powers between the test $TVE^{R}_{mn}$ and the tests $TV_{mn}$, $TVE_{mn}$, $TZ1_{mn}$, $TZ2_{mn}$ and $TC_{mn}$ are substantial. So, the test statistic $TVE^{R}_{mn}$ is most powerful test against the alternatives with the support $(-\infty,\infty)$.
\par
In groups III and IV, it is seen that the tests $TW^{R}_{mn}$ has the most power and the test $TVE_{mn}$ has the least power. Also, the difference of powers between $TV_{mn}$ and $TC_{mn}$ is small. But the difference of powers between the test $TW^{R}_{mn}$ and the tests $TV_{mn}$, $TVE_{mn}$, $TZ1_{mn}$, $TZ2_{mn}$  and $TC_{mn}$ are considerable. So, the test statistic $TW^{R}_{mn}$ is most powerful test against the alternatives with the support $(0,\infty)$.

   \item[]{\bf Case II: Test statistics in RSS context} \par
  {\textcolor{blue}{There already exist some goodness-of-fit tests in the RSS context. The most recent one appears to
   be Balakrishnan {\em et al.} (2013) \cite{bala2}, where they also considered the version of RSS Kolmogorov-Smirnov ($D^*$).
    For the  comparison, we evaluate corresponding powers of the test statistics $D^*$ and $S$ (
   which introduced by Balakrishnan {\em et al.} \cite{bala2}), $TVE^{R}_{mn}$ and $TW^{R}_{mn}$ by means of Monte Carlo 
   simulations for $8$ different alternative distributions. These alternatives were used by Balakrishnan {\em et al.} \cite{bala2}  
   as power comparison of several tests for normality. Table 13 contains the results of $10000$ simulations 
   (of samples size $25$, $30$, $45$, $50$, $60$, $75$, $90$ and $100$) per case to obtain the power of the 
   proposed tests and those of the competitor tests, at significance level $\alpha=0.05$. }}

\begin{sidewaystable}
\label{g3}
\caption{{\scriptsize Power comparisons based on $D^*$, $S$ in case II and $ECF$ in case III for several sample sizes under the alternatives with $\alpha=0.05$.}}
\vspace{0.1in}
\centering
\begin{tabular}{l|llll||llll||llll||llll}
\hline
\hline
\scriptsize Alt& \multicolumn{4}{c||}{\scriptsize Weibull(1.5)} & \multicolumn{4}{c||}{\scriptsize Weibull(2.5)} & \multicolumn{4}{c||}{\scriptsize chi(4)} & \multicolumn{4}{c}{\scriptsize Beta(2,2)}\\
\scriptsize $n$&\scriptsize $ECF$ &\scriptsize $D^*$ &\scriptsize S &\scriptsize $TW_{mn}$ &\scriptsize $ECF$ &\scriptsize $D^*$ &\scriptsize S &\scriptsize $TW_{mn}$ &\scriptsize $ECF$ &\scriptsize  $D^*$ &\scriptsize S &\scriptsize $TW_{mn}$  &\scriptsize $ECF$ &\scriptsize  $D^*$ &\scriptsize S &\scriptsize $TW_{mn}$ \\
\hline
 \scriptsize25&\scriptsize0.454 &\scriptsize  0.306 &\scriptsize 0.505 &\scriptsize \bf0.827 &\scriptsize 0.071   &\scriptsize 0.062 &\scriptsize 0.089 &\scriptsize \bf0.120 &\scriptsize  0.62  &\scriptsize 0.412 &\scriptsize 0.641 &\scriptsize\bf 0.926 &\scriptsize  0.040  &\scriptsize 0.043 &\scriptsize 0.045 &\scriptsize\bf 0.266 \\
 \scriptsize30&\scriptsize 0.533 &  \scriptsize0.362 & \scriptsize0.569 & \scriptsize\bf0.887 &\scriptsize  0.071  &\scriptsize 0.092 &\scriptsize 0.092 & \scriptsize \bf 0.133  &\scriptsize  0.688  &\scriptsize 0.499 &\scriptsize 0.713 & \scriptsize\bf0.963 &\scriptsize  0.042  & \scriptsize0.057 &\scriptsize 0.062 & \scriptsize\bf0.278 \\
 \scriptsize45&\scriptsize  0.786  &  \scriptsize0.471 & \scriptsize0.759 & \scriptsize\bf0.999 &\scriptsize  0.103  &\scriptsize 0.077 & \scriptsize0.106 & \scriptsize\bf0.154 &\scriptsize  0.865  & \scriptsize0.657 & \scriptsize0.895 &\scriptsize\bf 1.000 &\scriptsize 0.073   & \scriptsize0.075 &\scriptsize 0.092 &\scriptsize\bf 0.487 \\
 \scriptsize50&\scriptsize  0.808   &  \scriptsize0.565 & \scriptsize0.851 & \scriptsize\bf1.000 &\scriptsize 0.130   &\scriptsize 0.094 & \scriptsize0.126 & \scriptsize\bf0.164 &\scriptsize  0.903  & \scriptsize0.755 & \scriptsize0.937 &\scriptsize\bf 1.000 &\scriptsize  0.103  &\scriptsize 0.100 &\scriptsize 0.106 &\scriptsize \bf0.537 \\
 \scriptsize60&\scriptsize 0.869  &  \scriptsize0.627 & \scriptsize0.892 & \scriptsize\bf1.000 &\scriptsize  0.134  &\scriptsize 0.096 & \scriptsize0.137 & \scriptsize\bf0.263 &\scriptsize  0.946  & \scriptsize0.786 & \scriptsize0.966 &\scriptsize\bf 1.000 &\scriptsize  0.124  & \scriptsize0.084 &\scriptsize 0.147 &\scriptsize \bf0.733 \\
 \scriptsize75&\scriptsize  0.935    &  \scriptsize0.781 & \scriptsize0.961 & \scriptsize\bf1.000 &\scriptsize  0.180  &\scriptsize 0.126 & \scriptsize0.171 & \scriptsize\bf0.310 &\scriptsize  0.988  & \scriptsize0.908 & \scriptsize0.993 &\scriptsize \bf1.000 &\scriptsize  0.204  &\scriptsize 0.143 &\scriptsize 0.232 &\scriptsize \bf0.908 \\
 \scriptsize90&\scriptsize  0.976  &  \scriptsize0.810 & \scriptsize0.977 & \scriptsize\bf1.000 &\scriptsize 0.220   &\scriptsize 0.136 & \scriptsize0.218 & \scriptsize\bf0.438 &\scriptsize  0.989  & \scriptsize0.944 & \scriptsize0.995 &\scriptsize \bf1.000 &\scriptsize  0.266  &\scriptsize 0.144 &\scriptsize 0.297 &\scriptsize \bf0.986 \\
 \scriptsize100&\scriptsize 0.984  & \scriptsize0.890 & \scriptsize0.998 & \scriptsize\bf1.000 &\scriptsize  0.255  &\scriptsize 0.156 & \scriptsize0.225 &\scriptsize \bf0.404 &\scriptsize  0.997  & \scriptsize0.973 & \scriptsize0.999 & \scriptsize\bf1.000 &\scriptsize  0.318  &\scriptsize 0.169 &\scriptsize 0.358 & \scriptsize\bf0.992 \\
 \hline
\end{tabular}
\vspace{.1in}
\vspace{0.1in}
\centering
\begin{tabular}{l|llll||llll||llll||llll}
\hline
\scriptsize Alt& \multicolumn{4}{c||}{\scriptsize Beta(3,1)} & \multicolumn{4}{c||}{\scriptsize Beta(5,2)} & \multicolumn{4}{c||}{\scriptsize Uniform[0,1]} & \multicolumn{4}{c}{\scriptsize Logistic}\\
\scriptsize $n$&\scriptsize $ECF$&\scriptsize $D^*$ &\scriptsize S &\scriptsize $TW_{mn}$ &\scriptsize $ECF$&\scriptsize $D^*$ &\scriptsize S &\scriptsize $TW_{mn}$ &\scriptsize $ECF$&\scriptsize  $D^*$ &\scriptsize S &\scriptsize $TW_{mn}$  &\scriptsize $ECF$&\scriptsize  $D^*$ &\scriptsize S &\scriptsize $TVE_{mn} $ \\
\hline
 \scriptsize25&\scriptsize  0.460  &  \scriptsize0.372 &\scriptsize 0.595 &\scriptsize \bf0.995 &\scriptsize 0.170   &\scriptsize 0.146 &\scriptsize 0.261 &\scriptsize \bf0.416 &\scriptsize  0.144  &\scriptsize 0.099 & \scriptsize0.163 &\scriptsize\bf 0.896 &\scriptsize 0.122   &\scriptsize 0.095 &\scriptsize 0.122 &\scriptsize \bf0.141 \\
 \scriptsize30&\scriptsize 0.598   &  \scriptsize0.426 &\scriptsize 0.655 &\scriptsize \bf1.000 &\scriptsize 0.214   &\scriptsize 0.144 &\scriptsize 0.256 &\scriptsize \bf0.418 &\scriptsize  0.195  &\scriptsize 0.151 &\scriptsize 0.214 &\scriptsize\bf 0.956 &\scriptsize  0.143  & \scriptsize0.083 &\scriptsize 0.146 &\scriptsize\bf 0.162 \\
 \scriptsize45&\scriptsize  0.821  &  \scriptsize0.607 &\scriptsize 0.843 &\scriptsize \bf1.000 &\scriptsize 0.335   & \scriptsize0.236 &\scriptsize 0.398 &\scriptsize \bf0.784 &\scriptsize  0.419  &\scriptsize 0.216 &\scriptsize 0.431 &\scriptsize\bf 1.000 &\scriptsize  0.171  & \scriptsize0.097 &\scriptsize 0.177 & \scriptsize\bf0.220 \\
 \scriptsize50&\scriptsize  0.873  & \scriptsize 0.716 &\scriptsize 0.943 &\scriptsize\bf 1.000 &\scriptsize  0.398  & \scriptsize0.311 &\scriptsize 0.494 &\scriptsize \bf0.860 &\scriptsize  0.478  &\scriptsize 0.258 &\scriptsize 0.547 &\scriptsize \bf1.000 &\scriptsize  0.196  &\scriptsize 0.122 &\scriptsize 0.199 &\scriptsize\bf 0.257 \\
 \scriptsize60&\scriptsize 0.946   & \scriptsize 0.758 &\scriptsize 0.967 &\scriptsize \bf1.000 &\scriptsize 0.496   & \scriptsize0.324 &\scriptsize 0.574 &\scriptsize \bf0.960 &\scriptsize  0.618  &\scriptsize 0.313 &\scriptsize 0.654 & \scriptsize\bf1.000 &\scriptsize  0.199  & \scriptsize0.127 &\scriptsize 0.192 & \scriptsize\bf0.291 \\
 \scriptsize75&\scriptsize 0.976   & \scriptsize 0.920 &\scriptsize 0.992 &\scriptsize \bf1.000 &\scriptsize  0.607  &\scriptsize 0.430 &\scriptsize 0.716 &\scriptsize\bf 1.000 &\scriptsize  0.814  &\scriptsize 0.476 &\scriptsize 0.827 & \scriptsize\bf1.000 &\scriptsize  0.225  & \scriptsize0.163 & \scriptsize0.245 &\scriptsize \bf0.328 \\
 \scriptsize90&\scriptsize  0.996  &  \scriptsize0.937 &\scriptsize 0.995 &\scriptsize \bf1.000 &\scriptsize  0.733  & \scriptsize0.455 &\scriptsize 0.745 &\scriptsize\bf 1.000 &\scriptsize  0.882  &\scriptsize 0.534 &\scriptsize 0.920 &\scriptsize \bf1.000 &\scriptsize  0.251  & \scriptsize0.139 &\scriptsize 0.245 &\scriptsize \bf0.360 \\
\scriptsize 100&\scriptsize 0.998   & \scriptsize0.982 &\scriptsize 1.000 &\scriptsize \bf1.000 &\scriptsize 0.786   & \scriptsize0.547 &\scriptsize 0.847 &\scriptsize\bf 1.000 &\scriptsize  0.937  &\scriptsize 0.629 &\scriptsize 0.958 &\scriptsize \bf1.000&\scriptsize  0.261   & \scriptsize0.142 &\scriptsize 0.272 & \scriptsize\bf0.417 \\
 \hline
\end{tabular}
\end{sidewaystable}

  \item[]{\bf Case III: Some non-parametric test statistics} \par
  {\textcolor{blue} {For testing normality there are simple tests which are directly applied to the data at hand and do not use any transformation or RSS.
   As these methods are easier to apply ( in compare to our method where some parameters are to be fixed), we consider 
    power comparison with  such test statistics: ECF ( which introduced by Epps and Pulley \cite{epps}),
     Kolmogorov-Smirnov (KS), Anderson-Darling (AD) and Shapiro-Wilk  (SW). 
  We evaluated powers of the tests based on $ECF$, $TVE^{R}_{mn}$ and $TW^{R}_{mn}$ statistics, via Monte Carlo simulations, for the alternative distributions presented in Case II, that are reported in Table 13. Moreover, the power
  of test statistics $KS$, $SW$, $AD$, $TVE^{R}_{mn}$ and $TW^{R}_{mn}$, are computed for the alternative distributions presented in
   \cite{7test}, that  are reported  in Table 14. }}

\begin{sidewaystable}
\label{g3}
\caption{{\scriptsize Power comparisons based on $KS$, $AD$, $SW$ in case III for several sample sizes under the alternatives with $\alpha=0.05$.}}
\vspace{0.1in}
\centering
\begin{tabular}{l|llll||llll||llll||llll}
\hline
\hline
\scriptsize $n$ & \multicolumn{4}{c||}{\scriptsize $10$} & \multicolumn{4}{c||}{\scriptsize $20$} & \multicolumn{4}{c||}{\scriptsize $30$} & \multicolumn{4}{c}{\scriptsize $50$} \\
\scriptsize Alt. & \scriptsize $KS$ & \scriptsize $AD$ & \scriptsize $SW$ & \scriptsize $TW_{mn}^{R}$ & \scriptsize $KS$ & \scriptsize $AD$ & \scriptsize $SW$ & \scriptsize $TW_{mn}^{R}$ & \scriptsize $KS$ & \scriptsize $AD$ & \scriptsize $SW$ & \scriptsize $TW_{mn}^{R}$ & \scriptsize $KS$ & \scriptsize $AD$ & \scriptsize $SW$ & \scriptsize $TW_{mn}^{R}$ \\
\hline
\scriptsize Exp(1) & \scriptsize 0.301 & \scriptsize 0.416 & \scriptsize 0.442 & \scriptsize \bf0.537 & \scriptsize 0.586 & \scriptsize 0.773 & \scriptsize 0.836 & \scriptsize \bf0.999 & \scriptsize 0.784 & \scriptsize 0.934 & \scriptsize 0.968 & \scriptsize \bf1.000 & \scriptsize 0.963 & \scriptsize 0.997 & \scriptsize 0.999 & \scriptsize \bf1.000 \\
\scriptsize Gamma(2) & \scriptsize 0.210 & \scriptsize 0.225 & \scriptsize 0.239 & \scriptsize \bf0.250 & \scriptsize 0.326 & \scriptsize 0.467 & \scriptsize 0.532 & \scriptsize \bf0.744 & \scriptsize 0.472 & \scriptsize 0.660 & \scriptsize 0.749 & \scriptsize \bf0.960 & \scriptsize 0.701 & \scriptsize 0.890 & \scriptsize 0.949 & \scriptsize \bf1.000 \\
\scriptsize Gamma(0.5) & \scriptsize 0.672 & \scriptsize 0.703 & \scriptsize 0.735 & \scriptsize \bf0.916 & \scriptsize 0.879 & \scriptsize 0.970 & \scriptsize 0.984 & \scriptsize \bf1.000 & \scriptsize 0.982 & \scriptsize 0.998 & \scriptsize 0.999 & \scriptsize \bf1.000 & \scriptsize 0.999 & \scriptsize 1.000 & \scriptsize 1.000 & \scriptsize \bf1.000 \\
\scriptsize LN(0,1) & \scriptsize 0.463 & \scriptsize 0.578 & \scriptsize 0.603 & \scriptsize \bf0.738 & \scriptsize 0.778 & \scriptsize 0.904 & \scriptsize 0.932 & \scriptsize \bf1.000 & \scriptsize 0.935 & \scriptsize 0.984 & \scriptsize 0.991 & \scriptsize \bf1.000 & \scriptsize 0.995 & \scriptsize 0.999 & \scriptsize 0.999 & \scriptsize \bf1.000 \\
\scriptsize LN(0,2) & \scriptsize 0.826 & \scriptsize 0.909 & \scriptsize 0.920 & \scriptsize \bf0.995 & \scriptsize 0.991 & \scriptsize 0.998 & \scriptsize 0.999 & \scriptsize \bf1.000 & \scriptsize 0.999 & \scriptsize 1.000 & \scriptsize 1.000 & \scriptsize \bf1.000 & \scriptsize 1.000 & \scriptsize 1.000 & \scriptsize 1.000 & \scriptsize \bf1.000 \\
\scriptsize LN(0,0.5) & \scriptsize 0.182 & \scriptsize 0.233 & \scriptsize 0.245 & \scriptsize \bf0.250 & \scriptsize 0.337 & \scriptsize 0.463 & \scriptsize 0.517 & \scriptsize \bf0.649 & \scriptsize 0.492 & \scriptsize 0.656 & \scriptsize 0.726 & \scriptsize \bf0.890 & \scriptsize 0.715 & \scriptsize 0.871 & \scriptsize 0.924 & \scriptsize \bf0.999 \\
\scriptsize Weibull(0.5) & \scriptsize 0.758 & \scriptsize 0.875 & \scriptsize 0.894 & \scriptsize \bf0.990 & \scriptsize 0.981 & \scriptsize 0.997 & \scriptsize 0.999 & \scriptsize \bf1.000 & \scriptsize 0.999 & \scriptsize 1.000 & \scriptsize 1.000 & \scriptsize \bf1.000 & \scriptsize 1.000 & \scriptsize 1.000 & \scriptsize 1.000 & \scriptsize \bf1.000 \\
\scriptsize Weibull(2) & \scriptsize 0.074 & \scriptsize 0.083 & \scriptsize 0.084 & \scriptsize \bf0.090 & \scriptsize 0.103 & \scriptsize 0.132 & \scriptsize 0.156 & \scriptsize \bf0.176 & \scriptsize 0.132 & \scriptsize 0.184 & \scriptsize 0.232 & \scriptsize \bf0.271 & \scriptsize 0.210 & \scriptsize 0.307 & \scriptsize 0.416 & \scriptsize \bf0.701 \\
\scriptsize U(0,1) & \scriptsize 0.066 & \scriptsize 0.080 & \scriptsize 0.082 & \scriptsize \bf0.108 & \scriptsize 0.100 & \scriptsize 0.171 & \scriptsize 0.200 & \scriptsize \bf0.624 & \scriptsize 0.145 & \scriptsize 0.299 & \scriptsize 0.381 & \scriptsize \bf0.953 & \scriptsize 0.264 & \scriptsize 0.568 & \scriptsize 0.749 & \scriptsize \bf1.000 \\
\scriptsize Beta(2,2) & \scriptsize 0.046 & \scriptsize 0.046 & \scriptsize 0.042 & \scriptsize \bf0.066 & \scriptsize 0.053 & \scriptsize 0.058 & \scriptsize 0.053 & \scriptsize \bf0.158 & \scriptsize 0.060 & \scriptsize 0.080 & \scriptsize 0.080 & \scriptsize \bf0.264 & \scriptsize 0.083 & \scriptsize 0.127 & \scriptsize 0.153 & \scriptsize \bf0.654 \\
\scriptsize Beta(0.5,0.5) & \scriptsize 0.162 & \scriptsize 0.268 & \scriptsize 0.299 & \scriptsize \bf0.356 & \scriptsize 0.318 & \scriptsize 0.618 & \scriptsize 0.727 & \scriptsize \bf0.998 & \scriptsize 0.507 & \scriptsize 0.862 & \scriptsize 0.944 & \scriptsize \bf1.000 & \scriptsize 0.802 & \scriptsize 0.991 & \scriptsize 0.999 & \scriptsize \bf1.000 \\
\scriptsize Beta(3,0.5) & \scriptsize 0.418 & \scriptsize 0.576 & \scriptsize 0.609 & \scriptsize \bf0.819 & \scriptsize 0.746 & \scriptsize 0.913 & \scriptsize 0.948 & \scriptsize \bf1.000 & \scriptsize 0.934 & \scriptsize 0.991 & \scriptsize 0.997 & \scriptsize \bf1.000 & \scriptsize 0.998 & \scriptsize 1.000 & \scriptsize 1.000 & \scriptsize \bf1.000 \\
\scriptsize Beta(2,1) & \scriptsize 0.100 & \scriptsize 0.126 & \scriptsize 0.130 & \scriptsize \bf0.143 & \scriptsize 0.174 & \scriptsize 0.261 & \scriptsize 0.306 & \scriptsize \bf0.718 & \scriptsize 0.268 & \scriptsize 0.428 & \scriptsize 0.515 & \scriptsize \bf0.971 & \scriptsize 0.448 & \scriptsize 0.711 & \scriptsize 0.838 & \scriptsize \bf1.000 \\
 \hline
\end{tabular}
\vspace{.1in}
\vspace{0.1in}
\centering
\begin{tabular}{l|llll||llll||llll||llll}
\hline
\scriptsize $n$ & \multicolumn{4}{c||}{\scriptsize $10$} & \multicolumn{4}{c||}{\scriptsize $20$} & \multicolumn{4}{c||}{\scriptsize $30$} & \multicolumn{4}{c}{\scriptsize $50$} \\
\scriptsize Alt. & \scriptsize $KS$ & \scriptsize $AD$ & \scriptsize $SW$ & \scriptsize $TV_{mn}^{R}$ & \scriptsize $KS$ & \scriptsize $AD$ & \scriptsize $SW$ & \scriptsize $TV_{mn}^{R}$ & \scriptsize $KS$ & \scriptsize $AD$ & \scriptsize $SW$ & \scriptsize $TV_{mn}^{R}$ & \scriptsize $KS$ & \scriptsize $AD$ & \scriptsize $SW$ & \scriptsize $TV_{mn}^{R}$ \\
\hline
\scriptsize EV(0,1) & \scriptsize 0.121 & \scriptsize 0.147 & \scriptsize 0.153 & \scriptsize \bf0.219 & \scriptsize 0.203 & \scriptsize 0.273 & \scriptsize 0.313 & \scriptsize \bf0.350 & \scriptsize 0.290 & \scriptsize 0.402 & \scriptsize 0.469 & \scriptsize \bf0.473 & \scriptsize 0.440 & \scriptsize 0.596 & \scriptsize 0.686 & \scriptsize \bf0.727 \\
\scriptsize EV(0,2) & \scriptsize 0.121 & \scriptsize 0.144 & \scriptsize 0.150 & \scriptsize \bf0.152 & \scriptsize 0.203 & \scriptsize 0.276 & \scriptsize 0.315 & \scriptsize \bf0.396 & \scriptsize 0.289 & \scriptsize 0.399 & \scriptsize 0.464 & \scriptsize \bf0.474 & \scriptsize 0.437 & \scriptsize 0.579 & \scriptsize 0.690 & \scriptsize \bf0.721 \\
\scriptsize EV(0,0.5) & \scriptsize 0.118 & \scriptsize 0.147 & \scriptsize 0.154 & \scriptsize \bf0.173 & \scriptsize 0.203 & \scriptsize 0.275 & \scriptsize 0.314 & \scriptsize \bf0.335 & \scriptsize 0.288 & \scriptsize 0.400 & \scriptsize 0.464 & \scriptsize \bf0.474 & \scriptsize 0.440 & \scriptsize 0.598 & \scriptsize 0.690 & \scriptsize \bf0.708 \\
\scriptsize t(1) & \scriptsize 0.580 & \scriptsize 0.618 & \scriptsize 0.594 & \scriptsize \bf0.733 & \scriptsize 0.847 & \scriptsize 0.882 & \scriptsize 0.869 & \scriptsize \bf0.974 & \scriptsize 0.947 & \scriptsize 0.967 & \scriptsize 0.960 & \scriptsize \bf1.000 & \scriptsize 0.994 & \scriptsize 0.997 & \scriptsize 0.996 & \scriptsize \bf1.000 \\
\scriptsize t(3) & \scriptsize 0.164 & \scriptsize 0.190 & \scriptsize 0.187 & \scriptsize \bf0.207 & \scriptsize 0.260 & \scriptsize 0.327 & \scriptsize 0.340 & \scriptsize \bf0.446 & \scriptsize 0.345 & \scriptsize 0.436 & \scriptsize 0.460 & \scriptsize \bf0.554 & \scriptsize 0.484 & \scriptsize 0.599 & \scriptsize 0.632 & \scriptsize \bf0.839 \\
\scriptsize Logistic & \scriptsize 0.073 & \scriptsize 0.083 & \scriptsize 0.082 & \scriptsize \bf0.083 & \scriptsize 0.087 & \scriptsize 0.113 & \scriptsize 0.123 & \scriptsize \bf0.143 & \scriptsize 0.094 & \scriptsize 0.123 & \scriptsize 0.144 & \scriptsize \bf0.150 & \scriptsize 0.112 & \scriptsize 0.157 & \scriptsize 0.192 & \scriptsize \bf0.259 \\
\hline
\end{tabular}
\end{sidewaystable}

\end{description}

\section{Real Data}
\par
In this section we present a real examples to show the behavior of the test in real cases. For this example the parametric family used is the inverse Gaussian. This choice has been suggested in the literature since the example has been extensively analyzed.
\begin{example}{\em
The following real data represent active repair times (in hours) for an airborne communication transceiver.\\
 0.2, 0.3, 0.5, 0.5, 0.5, 0.5, 0.6, 0.6, 0.7, 0.7, 0.7, 0.8, 0.8, 1.0, 1.0, 1.0, 1.0, 1.1, 1.3, 1.5, 1.5,
1.5, 1.5, 2.0, 2.0, 2.2, 2.5, 3.0, 3.0, 3.3, 3.3, 4.0, 4.0, 4.5, 4.7, 5.0, 5.4, 5.4, 7.0, 7.5, 8.8, 9.0, 10.3,
22.0, 24.5.\\
The data were first analyzed by Von Alven \cite{von} who fitted successfully the Log normal distribution. Chhikara and Folks \cite{ha23} fitted the IG distribution and using the observed value of the Kolmogorov-Smirnov statistic, they concluded that the fit is good (KS test statistic $=0.07245267$). The same result is drawn by using the Anderson-Darling test (AD test statistic $=0.2392647$) and the Mudholkar independence characterization test (test statistic $=0.2026783$ \cite{ha24}). Finally, Lee {\em et al.} \cite{ha25} applied the IG distribution which was clearly accepted. The implementation of the proposed estimators can be used to investigate the above conclusion.
\par
Now, we apply the proposed goodness-of-fit test to the data set for testing $H_0: f_0(x)\sim N(\mu,\sigma^2)~~ \mbox{v.s}~~ H_1: f_1(x)\sim IG(\mu,\lambda)$ where $\mu,\sigma,\lambda$ are unknown. We use the data set and estimate: $\hat{\mu}=3.5644$, $\hat{\sigma}=4.7510$ and $\hat{\lambda}=2.0063$. The p-value of four test statistics can be earned in a simple method (see \cite{nigh} for details). {\textcolor{blue} { As $TW_{mn}^{R}$ is the most powerful against the alternatives with support $(0,\infty)$, we consider  $HW_{mn}^{R}$ as an entropy estimator. To choose the number of samples, $m$ and $w$, we use the Table 4. As an example for  $\mbox{MSE}\leq0.0214$, one could consider  $m=4$, $w=3$ and  $n\geq15$.
 To visualize the efficiency of this statistics we have evaluated the power of $TV_{mn}$, $TVE_{mn}$ and $TC_{mn}$ with $n=45$.  We find  that the power of the proposed test statistic, $TW_{mn}^{R}$, with $n=15$ is better than the power of the competitor test statistics with $n=45$.  Even though our method takes some more time, as for n=15 ( it takes 28 seconds  in compare to about 5 seconds for other methods), but the comparison of power which is 0.866 for the proposed test statistics $TW_{mn}^{R}$ to 0.5368, 0.1965, 4015 for other methods shows that the power of our method is significantly much higher.}}
 
{\textcolor{blue}{  We present an empirical procedure for estimating  powers of the entropy based test statistics, described in section 3. For this we provide some proper estimates
 of the corresponding parameters first.
  Then by simulating $M$ independent samples of size $n$ of the corresponding Inverse Gaussian distribution, we evaluate  corresponding test statistics.
 Finally we estimate powers of the test statistics for the comparison. This procedure can be described as:\\
\begin{description}
 \item[Step 1.] Fit a Normal distribution to the original data $x_1,\dots,x_n$ and obtain the estimated values $\hat{\mu}$, $\hat{\sigma}$.
 \item[Step 2.] Calculate the test statistic for the original data and denote by $TEST$.
 \item[Step 3.] Simulate $n$ independent sample $y_1,\dots,y_n$ of  $IG(\hat{\mu},\hat{\lambda})$ where $\hat{\lambda}=\frac{\hat{\mu}^3}{\hat{\sigma}^2}$.
 \item[Step 4.] Calculate corresponding test statistic for the simulated sample as $TEST^{1B}$.
 \item[Step 5.]  Repeat Step 3 and Step 4 for $M$ times and obtain $M$ test statistics as $TEST^{1B},\dots,TEST^{MB}.$
 \item[Step 6.] Evaluate the power by
 $\;\widehat{Power}=\frac{1}{M}\sum_{i=1}^{M}I_{(TEST^{iB}\leq TEST)},$
 where $I(.)$ is an indicator function taking value 1 when $TEST^{iB}\leq TEST$, and taking value 0 when $TEST^{iB}>TEST$.
\end{description}
By this algorithm, we estimate the power of different test statistics and report the results in Table 15.}}
\par
{\textcolor{blue} {Indeed, our results verify that, in all cases, the p-value of the different test statistics is much smaller than $0.05$ so that an IG distribution can provide a reasonable fitting to the data at $5\%$ level (Figure 2). Also, the power results show that the proposed test statistic, $TW^R_{mn}$, is most powerful test among the different test statistics, as expected. }}

\begin{table}
\caption{{\footnotesize P-value and power of different test statistics.}}
\centering
\vspace{0.1in}
\begin{tabular}{l|llll}
\hline
\hline
            & \scriptsize$TV_{mn}$ & \scriptsize$TVE_{mn}$ & \scriptsize$TC_{mn}$ & \scriptsize$TW_{mn}^R$\\
\hline
 \scriptsize p-value    & \scriptsize0.0011 & \scriptsize0.0000 & \scriptsize0.0002 & \scriptsize0.0028\\
 \scriptsize power  ($n=15$)    & \scriptsize0.5368 & \scriptsize0.1965 & \scriptsize0.4015 & \scriptsize0.8660\\
 \scriptsize power  ($n=45$)    & \scriptsize0.7175 & \scriptsize0.2905 & \scriptsize0.5019 & -\\
\hline
\end{tabular}
\end{table}

\label{fig3}\input{epsf}
\epsfxsize=3in \epsfysize=2in
 \begin{figure}
\centerline{\epsfxsize=4in \epsfysize=2.3in \epsffile{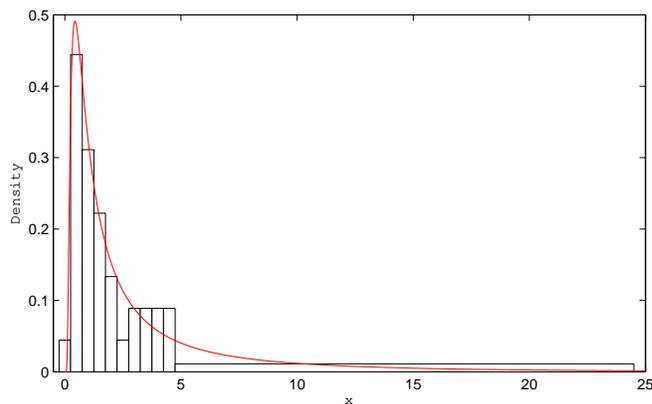}}
\vspace{-0.3in}
\caption{\footnotesize The histogram with an inverse gaussian fit for the active repair times. }
\end{figure}

}\end{example}
\section{Conclusion}
\par
In this paper, we first introduced two new entropy estimators of a continuous random variable by modifying the Van ES's \cite{31me} estimator and Wieczorkowski and Grzegorzewsky's \cite{33me} estimator in order to obtain two estimators with less bias and less RMSE. We considered some properties of $HVE^{R}_{mn}$ and $HW^{R}_{mn}$. 
Also, we compared our estimators with the entropy estimators proposed by Vasicek \cite{32me}, Van Es \cite{31me}, Ebrahimi {\em et al.} \cite{12me}, Correa \cite{8me}, Wieczorkowski and Grzegorewski \cite{33me} and Zamanzadeh and Arghami \cite{36me}. We observed that $HVE^{R}_{mn}$ behaves better than many of the most important entropy estimators and $HW^{R}_{mn}$ has generally least bias and RMSE among all entropy estimators for different distributions.
\par
We next introduced two new tests for normality based on the new estimators and discussed several tests of normality, using test statistics
\begin{itemize}
  \item $TV_{mn}$, $TVE_{mn}$, $TC_{mn}$, $TZ1_{mn}$, $TZ2_{mn}$, which are defined based on entropy estimators,
  \item $S$, $D^*$ which are defined in RSS context,
  \item $ECF$, $KS$, $AD$, $SW$, which are simple non-parametric test statistics.
\end{itemize}
 By comparing the power of these tests, using Monte Carlo computations, we found that the test statistics $TVE^{R}_{mn}$ and $TW^{R}_{mn}$ are most powerful against the alternatives with supports $(-\infty,\infty)$ and $(0,\infty)$, respectively.
\par
This work has the potential to be applied in the context of information theory and goodness-of-fit tests. This paper can elaborate further researches by considering other distributions besides the standard Normal distribution, such as Pareto, Log normal and Weibull distributions.


\bibliographystyle{plain}

\end{document}